\documentclass[11pt]{article}
\usepackage{amsfonts,amsmath,amssymb}
\usepackage{hyperref}
\newcommand{\status}{}

\newcommand{\tboc}{}

\newcommand{\file}{}

\newcommand{\detail}[1]{\par\noi{\bf [Proof detail\ }{#1}
\hfill{\bf ]}\par\noi\hspace{-4pt}}
\renewcommand{\detail}[1]{}

\newcommand{\dis}{\displaystyle}
\newcommand{\txt}{\textstyle}

\newcommand{\med}{\medskip}
\newcommand{\noi}{\noindent}

\newcommand{\halmos}{\rule{1ex}{1.4ex}}
\def \qed {\nopagebreak{\hspace*{\fill}$\halmos$\medskip}}

\newtheorem{theorem}{Theorem}[section]
\newtheorem{proposition}[theorem]{Proposition}
\newtheorem{corollary}[theorem]{Corollary}
\newtheorem{conjecture}[theorem]{Conjecture}
\newtheorem{lemma}[theorem]{Lemma}

\newtheorem{remark}[theorem]{Remark}

\newcommand{\bt}{\begin{theorem}}
\newcommand{\et}{\end{theorem}}
\newcommand{\bl}{\begin{lemma}}
\newcommand{\el}{\end{lemma}}
\newcommand{\bp}{\begin{proposition}}
\newcommand{\ep}{\end{proposition}}
\newcommand{\bcor}{\begin{corollary}}
\newcommand{\ecor}{\end{corollary}}
\newcommand{\br}{\begin{remark}\rm}
\newcommand{\er}{\end{remark}}
\newcommand{\bcon}{\begin{conjecture}}
\newcommand{\econ}{\end{conjecture}}

\newcommand{\be}{\begin{equation}}
\newcommand{\ee}{\end{equation}}
\newcommand{\ba}{\begin{array}}
\newcommand{\ea}{\end{array}}
\newcommand{\bc}{\be\begin{array}{r@{\,}c@{\,}l}}
\newcommand{\ec}{\end{array}\ee}

\newcommand{\ga}{\gamma}
\newcommand{\Ga}{\Gamma}
\newcommand{\de}{\delta}
\newcommand{\De}{\Delta}
\newcommand{\eps}{\varepsilon}
\newcommand{\la}{\lambda}
\newcommand{\La}{\Lambda}

\newcommand{\tet}{\theta}
\newcommand{\oo}{\omega}
\newcommand{\om}{\Omega}

\newcommand{\si}{\ensuremath{\sigma}}

\newcommand{\Ai}{{\cal A}}
\newcommand{\Bi}{{\cal B}}

\newcommand{\Di}{{\cal D}}
\newcommand{\Ei}{{\cal E}}
\newcommand{\Fi}{{\cal F}}
\newcommand{\Gi}{{\cal G}}

\newcommand{\Mi}{{\cal M}}

\newcommand{\Pc}{{\cal P}}
\newcommand{\Qi}{{\cal Q}}

\newcommand{\Si}{{\cal S}}

\newcommand{\Vi}{{\cal V}}

\newcommand{\R}{{\mathbb R}}
\newcommand{\N}{{\mathbb N}}
\newcommand{\Z}{{\mathbb Z}}

\newcommand{\E}{{\mathbb E}}
\renewcommand{\P}{{\mathbb P}}

\newcommand{\li}{\langle}
\newcommand{\re}{\rangle}

\newcommand{\volgt}{\ensuremath{\Rightarrow}}
\newcommand{\up}{\uparrow}
\newcommand{\down}{\downarrow}
\newcommand{\sub}{\subset}
\newcommand{\beh}{\backslash}

\newcommand{\symdif}{\!\vartriangle\!}

\newcommand{\asto}[1]{\underset{{#1}\to\infty}{\longrightarrow}}
\newcommand{\Asto}[1]{\underset{{#1}\to\infty}{\Longrightarrow}}
\newcommand{\astol}[1]{\underset{{#1}}{\longrightarrow}}

\newcommand{\ti}{\tilde}
\newcommand{\dgg}{\dagger}
\newcommand{\ov}{\overline}

\newcommand{\ffrac}[2]{{\textstyle\frac{{#1}}{{#2}}}}
\newcommand{\dif}[1]{\ffrac{\partial}{\partial{#1}}}

\newcommand{\di}{\mathrm{d}}
\newcommand{\half}{{[0,\infty)}}
\newcommand{\expo}{\mbox{\large\it e}}
\newcommand{\ex}[1]{\expo^{\,\textstyle{#1}}}

\newcommand{\var}{{\rm Var}}
\newcommand{\cov}{{\rm Cov}}

\newcommand{\expi}{\pi}
\newcommand{\poi}{\oo}

\begin{document}

%numbering formulas within sections
\makeatletter\@addtoreset{equation}{section}
\makeatother\def\theequation{\thesection.\arabic{equation}} 

%alternative layout for enumerate lists.
%\renewcommand{\labelenumi}{{\rm (\roman{enumi})}}

\title{\vspace{-3cm}The contact process
seen from a typical infected site}
\author{Jan M. Swart\vspace{6pt}\\
{\small Institute of Information}\\
{\small Theory and Automation}\\
{\small of the ASCR (\'UTIA)}\\
{\small Pod vod\'arenskou v\v e\v z\' i 4}\\
{\small 18208 Praha 8}\\
{\small Czech Republic}\\
{\small e-mail: swart@utia.cas.cz}\vspace{4pt}}
\date{{\small\file}June 16, 2008}
\maketitle\vspace{-.7cm}
\status

\begin{abstract}\noi
This paper studies contact processes on general countable groups. It
is shown that any such contact process has a well-defined exponential
growth rate, and this quantity is used to study the process. In
particular, it is proved that on any nonamenable group, the critical
contact process dies out.
\end{abstract}

\noi
{\it Acknowledgement.} Research supported by GA\v CR grant
201/06/1323 and the German Science Foundation. Part
of this work was carried out when the author was employed as a postdoc
at the university of T\"ubingen.

%60K35 Interacting random processes; statistical mechanics type models;
%      percolation theory
%82C22 Interacting particle systems
%82B43 Percolation

\tboc

\section{Introduction and main results}

\subsection{Introduction}

This paper studies contact processes whose underlying lattice is a
general countable group. There exists a small body of literature about
contact processes on general lattices, but several basic questions
have been answered only on specific lattices. In particular, a lot is
known about the process on the $d$-dimensional integer lattice $\Z^d$,
and on regular trees. (See \cite{Lig99} as a general reference for
contact processes on $\Z^d$, trees, and other lattices.)

It turns out that the contact process on regular trees behaves quite
differently from the contact process on $\Z^d$. For the process on $\Z^d$, it
is known that there is a critical infection rate $0<\la_{\rm c}<\infty$ such
that for $\la\leq\la_{\rm c}$, the process dies out, while for $\la>\la_{\rm
  c}$, the process survives with positive probability, and complete
convergence holds. On the other hand, on trees, there are two critical values
$0<\la_{\rm c}<\la'_{\rm c}<\infty$ such that in the intermediate regime
$\la_{\rm c}<\la\leq\la'_{\rm c}$, the process survives, but complete
convergence does not hold. The situation is quite similar to the situation for
(unoriented) percolation on general transitive lattices, where it is known
that one has uniqueness of the infinite cluster whenever the lattice is
amenable, while it is conjectured (and proved in several special cases) that
on any nonamenable lattice there exists an intermediate parameter regime where
there are infinitely many infinite clusters.

While a lot is known nowadays about percolation on general transitive
graphs, the same cannot be said for the contact process. In
particular, it is not known what is the essential difference between
$\Z^d$ and trees that causes the observed difference in behavior on
these lattices. A natural guess is that the essential feature is
amenability ($\Z^d$ being amenable, while trees are not). However,
as we will see shortly, there are reasons to doubt this.

In the present paper, we study contact processes on general countable groups
by means of their exponential growth rate. A simple subadditivity argument
shows that the expected number of infected sites of a contact process on a
transitive lattice, started with finitely many infected sites, grows at a
well-defined exponential rate (independent of the initial state). On $\Z^d$,
it is known that this exponential growth rate is negative for $\la<\la_{\rm
  c}$ (see \cite{BG91} or \cite[Thm~I.2.48]{Lig99}), and zero for any
$\la\geq\la_{\rm c}$. Indeed, it is easy to see (and prove) that on $\Z^d$
there is simply not enough space for a contact process to grow exponentially
fast (with positive exponent). On the other hand, one of our main results in
this paper is that if a contact process survives on a nonamenable group, then
its exponential growth rate must be strictly positive. This result is known
for trees; our proof in the case of general nonamenable groups is quite
different from the known proof for trees, however. The main idea of our proof
is to relate the exponential growth rate of a contact process to the
configuration seen by a typical infected site at a typical late time.

Intuition says that a contact process that survives with a positive
exponential growth rate behaves very much like a perturbed branching
process. On the other hand, contact processes that survive but have a zero
exponential growth rate are different. We do not know if (non)amenability is
the essential feature here. It is known that there exist exponentially growing
groups that are amenable. (A well-known example is the lamplighter
group). Although we do not prove it here, it seems plausible that a contact
process on such a group, if it survives, must have a positive exponential
growth rate. Thus, contact processes on such amenable groups might in some
respects show behavior that is more similar to processes on trees than on
$\Z^d$. It should be noted that (at least) on non-homogeneous lattices, the
situation is even more complex. In particular, Pemantle and Stacey \cite{PS01}
have given examples of non-homogeneous trees of uniformly exponential growth
and bounded degree, on which the critical values related to survival and
complete convergence of the contact process coincide.

Part of the present work appeared before as Chapter~4 of the author's
habilitation thesis \cite{Swa07}. In particular,
Proposition~\ref{P:Campconv} below is Theorem~4.3~(a) in \cite{Swa07}.

\subsection{Set-up}

We will study contact processes whose underlying lattice is a general
countable group. From the point of view of studying general
transitive lattices, this is not quite as general as one might
wish; in particular, such lattices are always unimodular. Assuming
that the lattice is a group will simplify our proofs, however, so as a
first step it seems reasonable.

Our set-up is as follows. We let $\La$ be a finite or countably
infinite group, which we refer to as the {\em lattice}, with group
action $(i,j)\mapsto ij$ and unit element $0$, also referred to as the
origin. Each site $i\in\La$ can be in one of two states: healthy or
infected. Infected sites become healthy with {\em recovery rate}
$\de\geq 0$. An infected site $i$ infects another site $j$ with {\em
infection rate} $a(i,j)\geq 0$. We assume that the infection rates are
invariant with respect to the left action of the group and summable:
\be\ba{rl}\label{assum}
{\rm (i)}&a(i,j)=a(ki,kj)\qquad\qquad(i,j,k\in\La),\\[5pt]
{\rm (ii)}&\dis|a|:=\sum_ia(0,i)<\infty,\\[5pt]
\ec
Here we adopt the convention that sums over $i,j,k$ always run over
$\La$, unless stated otherwise. Note that we do {\em not} assume that
$a(i,j)=a(j,i)$, i.e., our contact processes are in general asymmetric.

Let $\eta_t$ be the set of all infected sites at time $t\geq 0$. Then
$\eta=(\eta_t)_{t\geq 0}$ is a Markov process in the space
$\Pc(\La):=\{A:A\sub\La\}$ of all subsets of $\La$, called the contact
process on $\La$ with infection rates $a=(a(i,j))_{i,j\in\La}$ and
recovery rate $\de$, or shortly the {\em $(\La,a,\de)$-contact
process}. We equip $\Pc(\La)\cong\{0,1\}^\La$ with the product
topology and the associated Borel-\si-field $\Bi(\Pc(\La))$, and let
$\Pc_{\rm fin}(\La):=\{A\sub\La:|A|<\infty\}$ denote the subspace of
finite subsets of $\La$. Under the assumptions (\ref{assum}), $\eta$
is a well-defined Feller process with cadlag sample paths in the
compact state space $\Pc(\La)$, and $\eta_0\in\Pc_{\rm fin}(\La)$
implies $\eta_t\in\Pc_{\rm fin}(\La)$ for all $t\geq 0$ a.s.

Note that we have not assumed any additional structure on $\La$,
except for the group structure. In particular, we have not assumed any
sort of `nearest neighbor' structure. This may be obtained in the
following special case. Assume that $\La$ is finitely generated and
that $\De$ is a finite, symmetric (with respect to taking inverses),
generating set for $\La$. Then the (left) Cayley graph
$\Gi=\Gi(\La,\De)$ associated with $\La$ and $\De$ is the graph with
vertex set $\Vi(\Gi):=\La$ and edges
$\Ei(\Gi):=\{\{i,j\}:i^{-1}j\in\De\}$. Examples of Cayley graphs are
$\Z^d$ and regular trees. (In the case of trees, there are several
possible choices for the group structure.) Setting
$a(i,j):=\la1_{\{i^{-1}j\in\De\}}$, with $\la>0$, and choosing
$\de\geq 0$, then defines a nearest-neighbor contact process on the
Cayley graph $\Gi(\La,\De)$. In this case, $\la$ is simply referred to
as `the' infection rate. If $\de>0$, then by rescaling time we may set
$\de=1$, so it is customary to assume that $\de=1$. If $\de=0$, then
$\eta$ is a special case of first-passage percolation (see
\cite{Kes86}).

Returning to our more general set-up, we make the following
observation, which is the basis of our analysis. Below, we use the
notation $\eta^A_t$ to denote the $(\La,a,\de)$-contact process
started at time zero in $\eta^A_0=A$, evaluated at time $t\geq 0$.
\bl{\bf(Exponential growth rate)}\label{L:exprate}
Let $\eta$ be a $(\La,a,\de)$-contact process. Then there exists a
constant $r=r(\La,a,\de)$ with $-\de\leq r\leq|a|-\de$ such that
\be\label{expr}
\lim_{t\to\infty}\,\ffrac{1}{t}\log\E\big[|\eta^A_t|\big]=r
\qquad(\emptyset\neq A\in\Pc_{\rm fin}(\La)).
\ee
\el
We call $r=r(\La,a,\de)$ the {\em exponential growth rate} of the
$(\La,a,\de)$-contact process. We note that $r$ has been defined
before in the specific context of nearest-neighbor processes on
regular trees. Indeed, $r=\log\phi(1)$, where $\phi(\rho)$ is the
function defined in \cite[formula~(I.4.23)]{Lig99}.

\subsection{The exponential growth rate}

In this section, we investigate the exponential growth rate
$r(\La,a,\de)$ of a contact process defined in Lemma~\ref{L:exprate}.

We start by recalling a few basic facts and definitions concerning
groups. As before, let $\La$ be a finite or countably infinite group.
For $i\in\La$ and $A,B\sub\La$ we put $AB:=\{ij:i\in A,\ j\in B\}$,
$iA:=\{i\}A$, $Ai:=A\{i\}$, $A^{-1}:=\{i^{-1}:i\in A\}$, $A^0:=\{0\}$,
$A^n:=AA^{n-1}$ $(n\geq 1)$, and $A^{-n}:=(A^{-1})^n=(A^n)^{-1}$. We
write $A\symdif B:=(A\beh B)\cup(B\beh A)$ for the symmetric difference
of $A$ and $B$ and let $|A|$ denote the cardinality of $A$.

By definition, we say that $\La$ is {\em amenable} if
\be\ba{l}\label{amen}
\mbox{For every finite nonempty $\De\sub\La$ and $\eps>0$,
there exists a finite}\\
\mbox{nonempty $A\sub\La$ such that $|(A\De)\symdif A|\leq\eps|A|$.}
\ec
If $\La$ is finitely generated, then it suffices to check (\ref{amen})
for one finite symmetric generating set $\De$. In this case, $(A\De)\symdif
A$ is the set of all $i\not\in A$ for which there exists a $j\in A$
such that $i$ and $j$ are connected by an edge in the Cayley graph
$\Gi(\La,\De)$. Thus, we may describe (\ref{amen}) by saying that it
is possible to find nonempty sets $A$ whose surface is small compared
to their volume. For example, $\Z^d$ is amenable, but regular trees
are not.

If $\La$ is a finitely generated group and $\De$ is a finite symmetric
generating set, then we let $|i|$ denote the usual graph distance of
$i$ to the origin in the Cayley graph $\Gi(\La,\De)$, i.e.,
$|i|:=\min\{n:i\in\De^n\}$.  The norm $|\cdot|$ depends on the choice
of $\De$, but any two norms associated with different finite symmetric
generating sets are equivalent. It follows from subadditivity that the
limit $\lim_{n\to\infty}\frac{1}{n}\log|\{i\in\La:|i|\leq n\}|$ exists;
one says that the group $\La$ has exponential (resp.\ subexponential)
growth if this limit is positive (resp.\ zero). Note that since norms
associated with different finite symmetric generating sets are equivalent,
having (sub)exponential growth is a property of the group $\La$ only
and does not depend on the choice of~$\De$.

Subexponential growth implies amenabilty, but the converse is not
true: as already mentioned, the lamplighter group is an
amenable group with exponential growth. See \cite[Section~5]{MW89} for
general facts about amenability and subexponential growth, and
\cite{LPP96} for a nice exposition of the lamplighter group.

We also need a few definitions concerning contact processes. If
$a=(a(i,j))_{i,j\in\La}$ are infection rates satisfying (\ref{assum}),
then we define {\em reversed infection rates} $a^\dgg$ by
$a^\dgg(i,j):=a(j,i)$ $(i,j\in\La)$. We say that the
$(\La,a,\de)$-contact process {\em survives} if
\be\label{surv}
\P\big[\eta^A_t\neq\emptyset\ \forall t\geq 0\big]>0
\ee
for some, and hence for all $\emptyset\neq A\in\Pc_{\rm fin}(\La)$.
Using the standard coupling, it is easy to see that if $\de<\de'$ and
the $(\La,a,\de')$-contact process survives, then the
$(\La,a,\de)$-contact process survives. We let
\be
\de_{\rm c}=\de_{\rm c}(\La,a):=\sup\big\{\de\geq 0:\mbox{ the
$(\La,a,\de)$-contact process survives}\big\}
\ee
denote the {\em critical recovery rate}. By comparison with a critical
branching process, it is not hard to see that $\de_{\rm c}\leq|a|$. Although
we do not need this in what follows, we note that if $\La$ is finitely
generated and the infection rates $a$ are irreducible, then one may use
comparison with a one-dimensional nearest-neighbor contact process to show
that $0<\de_{\rm c}$ (see \cite[Lemma~4.18]{Swa07}).

Here, we say that infection rates $a$ on $\La$ are {\em irreducible} if
\be\label{ir1}
\bigcup_{n\geq 0}(A\cup A^{-1})^n=\La\quad\mbox{where }A:=\{i\in\La:a(0,i)>0\}.
\ee
At some point, we will need an assumption that is a bit stronger than
this. More precisely, we will occasionally use the following
assumption (see Lemma~\ref{L:evdom} below):
\be\label{irr}
\bigcup_{n\geq 0,\ m\geq 0}\!\!A^{-n}A^m=\La
=\bigcup_{n\geq 0,\ m\geq 0}\!\!A^nA^{-m},
\quad\mbox{where }A:=\{i\in\La:a(0,i)>0\}.
\ee
Note that this says that for any two sites $i,j$ there exists
a site $k$ from which both $i$ and $j$ can be infected, and a site $k'$
that can be infected both from $i$ and from $j$.

With these definitions, we are ready to formulate our main result.
\bt{\bf(Properties of the exponential growth rate)}\label{T:main}
Let $\La$ be a finite or countably infinite group, let $a=(a(i,j))_{i,j\in\La}$
be infection rates satisfying (\ref{assum}), and $\de\geq 0$. Let
$r=r(\La,a,\de)$ be the exponential growth rate of the
$(\La,a,\de)$-contact process, defined in (\ref{expr}). Then:
\begin{itemize}
\item[\rm(a)] $r(\La,a,\de)=r(\La,a^\dgg,\de)$

\item[\rm(b)] The function $\de\to r(\La,a,\de)$ is nonincreasing and
Lipschitz continuous on $\half$, with Lipschitz constant 1.

\item[\rm(c)] If the $(\La,a,\de)$-contact process survives, then $r\geq 0$.

\item[\rm(d)] If $r>0$, then the $(\La,a,\de)$-contact process survives.

\item[\rm(e)] If $\La$ is finitely generated and has subexponential
growth, and the infection rates satisfy
$\sum_ia(0,i)e^{\eps|i|}<\infty$ for some $\eps>0$, then $r\leq 0$.

\item[\rm(f)] If $\La$ is nonamenable, the $(\La,a,\de)$-contact process
survives, and the infection rates satisfy the irreducibility condition
(\ref{irr}), then $r>0$.
\end{itemize}
\et
Parts~(a), (b), and (c) of this theorem are easy. Part~(d) follows from a
variance calculation, while (e) is proved by some simple large deviation
estimates. The proof of part~(f) is rather involved. For trees, the statement
is a known consequence of \cite[Prop.~I.4.27~(b)]{Lig99}. Our proof for
general nonamenable groups is quite different from the methods used there.
The basic idea is as follows. If the exponential growth rate of a contact
process is zero, then this means that for the process started with one
infected site, a `typical' infected site at a `typical' late time produces no
net offspring, i.e., the mean number of sites it infects per unit of time is
just enough to balance the probability that the site itself recovers. We prove
that this implies that the local configuration as seen from this `typical'
site is distributed as the upper invariant measure, assuming that the latter
is nontrivial. If $\La$ is nonamenable, this leads to a contradiction, since
for any finite collection of particles on a nonamenable lattice, a positive
fraction of the particles must lie on the `outer boundary' of the collection,
hence must see something different from the upper invariant measure.

Parts~(b), (d), and (f) of Theorem~\ref{T:main} yield the following corollary.
\bcor{\bf(The critical contact process on a nonamenable lattice dies out)}
\label{C:crit}
If $\La$ is nonamenable and the infection rates satisfy the
irreducibility condition (\ref{irr}), then $\de_{\rm c}=\de_{\rm
c}(\La,a)>0$ and the $(\La,a,\de_{\rm c})$-contact process dies out.
\ecor
For nearest-neighbor contact processes on regular trees, this result is known,
see \cite{MSZ94} or \cite[Proposition~I.4.39]{Lig99}. Like our proof, the
proof there is based on showing that a zero exponential growth rate implies
extinction (although they use quite different techniques to establish
this). The analogue of Corollary~\ref{C:crit} for (unoriented) percolation
says that there are no infinite clusters at criticality on any nonamenable
lattice. This has been proved in \cite{BLPS99}; it seems that the techniques
used there have little in common with the ones used in the present paper.  The
problem of showing that critical percolation on $\Z^d$ has no infinite cluster
is still open in dimensions $3\leq d\leq 18$. Barsky, Grimmett, and Newman
\cite{BGN91} have shown, however, that at criticality there are no infinite
clusters in the half-space $\Z_+\times\Z^{d-1}$, and, using similar techniques,
Bezuidenhout and Grimmett \cite{BG90} have proved that the critical contact
process on $\Z^d$ dies out.

\subsection{The process seen from a typical site}

There is an intimate relation between the survival probability of a
contact process and its upper invariant law. Similarly, there is a
relation between the exponential growth rate and certain infinite
measures on the space of nonempty subsets of $\La$, which we explain
now.

Let $\Pc_+(\La):=\{A\sub\La:A\neq\emptyset\}$ denote the set of all nonempty
subsets of $\La$. Note that $\Pc_+(\La)$ is a locally compact space in the
induced topology from $\Pc(\La)$. We say that a measure $\mu$ on $\Pc(\La)$ or
$\Pc_+(\La)$ is (spatially) {\em homogeneous} if it is invariant under the
left action of the group, i.e., if $\mu(\Ai)=\mu(i\Ai)$ for each $i\in\La$ and
$\Ai\in\Bi(\Pc(\La))$, where we define $i\Ai:=\{iA:A\in\Ai\}$. We say that a
measure $\mu$ is an {\em eigenmeasure} of the $(\La,a,\de)$-contact process if
$\mu$ is a nonzero, locally finite measure on $\Pc_+(\La)$, and there exists a
constant $\la\in\R$ such that
\be\label{eigen}
\int\mu(\di A)\P[\eta^A_t\in\cdot\,]\big|_{\Pc_+(\La)}
=e^{\la t}\mu\qquad(t\geq 0),
\ee
where $|_{\Pc_+(\La)}$ denotes restriction (of a measure) to $\Pc_+(\La)$. We
call $\la$ the associated {\em eigenvalue}. As a motivation for this
terminology, we observe that if $G$ is the generator of the
$(\La,a,\de)$-contact process, then formally $G^\ast\mu=\la\mu$. Note that if
$\la=0$ and $\mu$ is concentrated on the infinite subsets of $\La$, then the
measure on the left-hand side of (\ref{eigen}) is concentrated on
$\Pc_+(\La)$, hence in this case (\ref{eigen}) just says that $\mu$ is an
invariant measure (though not necessarily a probability measure) for the
$(\La,a,\de)$-contact process.

\bp{\bf(Exponential growth rate and eigenmeasures)}\label{P:eigen}
For each $(\La,a,\de)$-con\-tact process, the set
\be\ba{r@{\,}l}
\dis\Ei(\La,a,\de):=\big\{\la\in\R:&\dis\mbox{there exists a homogeneous
eigenmeasure}\\[5pt]
&\dis\mbox{of the $(\La,a,\de)$-contact process with eigenvalue }\la\big\}
\ec
is a nonempty compact subset of $\R$, and $r(\La,a,\de)=\max\Ei(\La,a,\de)$.
\ep
In particular, Proposition~\ref{P:eigen} implies that each
$(\La,a,\de)$-contact process has a homogeneous eigenmeasure
with eigenvalue $r(\La,a,\de)$. It seems natural to conjecture that this
eigenmeasure is always unique and the long-time limit of the (suitably
rescaled) law of the process started with one infected site, distributed
according to counting measure on $\La$. If we condition such an eigenmeasure
on the origin being infected, then we can view the resulting probability
measure as describing the contact process as seen from a typical infected
site, at late times. (Compare Corollary~\ref{C:conv} and Lemma~\ref{L:Camp}
below.)

Recall that the upper invariant measure $\ov\nu$ of a contact process is the
long-time limit law of the process started with all sites infected. It follows
{f}rom duality (and spatial ergodicity of the graphical representation) that the
upper invariant measure of the $(\La,a,\de)$-contact process is nontrivial
(i.e., gives zero probability to the empty set) if and only if the
$(\La,a^\dgg,\de)$-contact process survives. The next result is an important
ingredient in the proof of Theorem~\ref{T:main}~(f).
\bt{\bf(Eigenmeasures with eigenvalue zero)}\label{T:zero}
Assume that the infection rates satisfy the irreducibility condition
(\ref{irr}). If the upper invariant measure $\ov\nu$ of the
$(\La,a,\de)$-contact process is nontrivial, then any homogeneous
eigenmeasure $\mu$ with eigenvalue zero satisfies
$\mu=c\ov\nu$ for some $c>0$.
\et
We prove Theorem~\ref{T:zero} by extending well-known techniques for showing
that $\ov\nu$ is the only nontrivial homogeneous invariant probability measure
of a contact process. If $\ov\nu$ is trivial, then there may exist
homogeneous eigenmeasures with eigenvalue zero, which in this
case, obviously, are not a multiple of $\ov\nu$. Indeed, if $a$ is symmetric
(i.e., $a=a^\dgg$) and the critical process dies out (as we know to be the
case on $\Z^d$ or on any nonamenable group), then at criticality $\ov\nu$ is
trivial, while by Theorem~\ref{T:main} (b), (c), and (d), the exponential growth
rate is zero, hence by Proposition~\ref{P:eigen}, there exists a homogeneous
eigenmeasure with eigenvalue zero.

\subsection{Discussion, open problems, and outline}

The work in this paper started from the question whether it is
possible to prove something like `uniqueness of the infinite cluster'
in the context of oriented percolation or the (very similar) graphical
representation of the contact process. This question is still very
much open. See Grimmett and Hiemer \cite{GH02} for a weak statement
that is proved only on $\Z^d$ and Wu and Zhang \cite[Thm 1.4]{WZ06} or
\cite[Lemma~4.5]{Swa07} for a stronger statement that is proved only
in the nearest-neighbor, one-dimensional case.

Whether the methods in the present paper can shed some light on this
question I do not know. I have tried to prove the weak statement of
Grimmett and Hiemer assuming (only) subexponential growth, but ran
into the problem that I would need to replace a size-biased law by a
law conditioned on survival, which I do not know how to do (see
\cite[Prop~4.4]{Swa07}).

In fact, although this is not obvious from the presentation above,
size-biased laws and Campbell measures, well-known objects from
branching theory, are closely related to the eigenmeasures introduced
above. (For this connection, see Section~\ref{S:Palm} below.) An
interesting feature of the (potentially infinite) eigenmeasures is that
they allow one to use some of the simplifications that come from
spatial homogeneity while studying processes started in finite initial
states.

There are lots of open problems concerning contact processes on
general transitive lattices, so we mention just a few.
\begin{enumerate}
\item Prove that the $(\La,a,\de)$-contact process has a unique
homogeneous eigenmeasure with eigenvalue
$r(\La,a,\de)$, which is the long-time limit law of the process
started with one infected site distributed according to counting
measure on $\La$.

\item Prove that $\dif{\de}r(\La,a,\de)<0$ on $\{\de:r(\La,a,\de)\neq
  0\}$. Prove the same statement for all $\de$ if $\La$ is nonamenable. Adapt
  the known proof for $\Z^d$ (see \cite{BG91} and \cite[Thm~I.2.48]{Lig99})
  that $r(\La,a,\de)<0$ for all $\de>\de_{\rm c}$, to general lattices.

\item Prove that $\de_{\rm c}>0$ for some $(\La,a,\de)$-contact
process on a group $\La$ that is not finitely generated, e.g.\ the
hierarchical group.

\item Study contact processes on transitive lattices $\La$ that are
not groups. In this context, if $\La$ is not unimodular, it is not hard to
find examples where a $(\La,a,\de)$-contact process survives but its dual
$(\La,a^\dgg,\de)$-contact process dies out. It is an open problem to prove
this cannot happen in the unimodular case.

\item Prove some version of uniqueness of the infinite cluster
assuming that the exponential growth rate is zero.

\item Prove (or disprove) that $r(\La,a,\de)>0$ whenever the
  $(\La,a,\de)$-contact process survives and $r(\La,a,0)>0$.
\end{enumerate}
The outline of the rest of the paper is as follows. In
Section~\ref{S:const}, we introduce some basic tools, such as the
graphical representation and a martingale problem. In
Section~\ref{S:exp}, we prove Theorem~\ref{T:main}~(a)--(c),
Proposition~\ref{P:eigen}, and Theorem~\ref{T:zero}. These results are
then used in Section~\ref{S:proof} to prove
Theorem~\ref{T:main}~(d)--(f) and Corollary~\ref{C:crit}.\med

\noi
{\bf\large Acknowledgements} The author thanks Geoffrey Grimmmett, Olle
H\"aggstr\"om, Russel Lyons, Amos Nevo, Yuval Peres, and Roberto Schonmann for
useful email conversations about the contact process, oriented percolation,
and amenability. In particular, the proof of formula (\ref{unizero}) is due to
Yuval Peres. The author thanks the referee for many useful suggestions.
%for improving the paper.

\section{Construction and basic properties}\label{S:const}

\subsection{Graphical representation}\label{S:graph}

We will, of course, use the graphical representation
of the contact process. Let $\La\times\R:=\{(i,t):i\in\La,\ t\in\R\}$
and $\La\times\La\times\R:=\{(i,j,t):i,j\in\La,\ t\in\R\}$, where $t$
is the time coordinate. Let $\poi=(\poi^{\rm r},\poi^{\rm i})$ be a
pair of independent, locally finite random subsets of $\La\times\R$
and $\La\times\La\times\R$, respectively, produced by Poisson point
processes with intensity $\de$ and $a(i,j)$, respectively. We
visualize this by plotting $\La$ horizontally and $\R$ vertically,
marking points $(i,s)\in\poi^{\rm r}$ with a recovery symbol $\ast$,
and drawing an infection arrow from $(i,t)$ to $(j,t)$ for each
$(i,j,t)\in\poi^{\rm i}$. For $C,D\sub\La\times\R$, say that there is
a {\em path} from $C$ to $D$, denoted by $C\leadsto D$, if there exist
$n\geq 0$, $i_0,\ldots,i_n\in\La$, and $t_0\leq\cdots\leq t_{n+1}$
with $(i_0,t_0)\in C$ and $(i_n,t_{n+1})\in D$, such that
$(\{i_k\}\times[t_k,t_{k+1}])\cap\poi^{\rm r}=\emptyset$ for all
$k=0,\ldots,n$ and $(i_{k-1},i_k,t_k)\in\poi^{\rm i}$ for all
$k=1,\ldots,n$. Thus, a path must walk upwards in time, may follow
arrows, and must avoid recoveries. For $C\sub\La\times\R$, we write
$C\leadsto\infty$ if there is an infinite path with times
$t_k\up\infty$ starting in $C$. We define $-\infty\leadsto C$
analogously. Instead of $\{(i,s)\}\leadsto\:$ and
$\:\leadsto\{(j,t)\}$, simply write $(i,s)\leadsto\:$ and
$\:\leadsto(j,t)$.

For given $A\in\Pc(\La)$ and $t_0\in\R$, put
\be\label{eta}
\eta^{A\times\{t_0\}}_t:=\{i\in\La:A\times\{t_0\}\leadsto(i,t_0+t)\}
\qquad(t\geq 0).
\ee
Then $\eta^{A\times\{t_0\}}=(\eta^{A\times\{t_0\}}_t)_{t\geq 0}$ is a
$(\La,a,\de)$-contact process started in $\eta^{A\times\{t_0\}}_0=A$.

In analogy with (\ref{eta}), put
\be
\eta^{\dgg\,A\times\{t_0\}}_t:=\{i\in\La:(i,t_0-t)\leadsto A\times\{t_0\}\}
\qquad(t\geq 0).
\ee
Then $\eta^{\dgg\,A\times\{t_0\}}=(\eta^{\dgg\,A\times\{t_0\}}_t)_{t\geq 0}$
is a $(\La,a^\dgg,\de)$-contact process started in
$\eta^{\dgg\,A\times\{t_0\}}_0=A$. Since for any $s\leq t$ and
$A,B\in\Pc(\La)$, the event
\be
\big\{\eta^{A\times\{s\}}_{u-s}\cap\eta^{\dgg\,B\times\{t\}}_{t-u}
=\emptyset\big\}=\big\{A\times\{s\}\not\leadsto B\times\{t\}\big\}
\ee
does not depend on $u\in[s,t]$, it follows (by taking $s=0$ and $u=0,t$) that
the $(\La,a,\de)$-contact process and the $(\La,a^\dgg,\de)$-contact process
are dual in the sense that
\be\label{dual}
\P[\eta^A_t\cap B=\emptyset]=\P[A\cap\eta^{\dgg\,B}_t=\emptyset]
\qquad(A,B\in\Pc(\La),\ t\geq 0).
\ee
Here, for brevity, we write $\eta^A_t:=\eta^{A\times\{0\}}_t$ and
$\eta^{\dgg\,B}_t:=\eta^{\dgg\,B\times\{0\}}_t$.

It is not hard to see that $|a|:=\sum_ia(0,i)=\sum_ia(i,0)$ and
\be\label{fiex}
\E\big[|\eta^A_t|\big]\leq|A|e^{|a|t}\quad\mbox{and}\quad
\E\big[|\eta^{\dgg\,A}_t|\big]\leq|A|e^{|a|t}\qquad(t\geq 0,\ A\sub\La).
\ee
In particular, both the $(\La,a,\de)$-contact process and the
$(\La,a^\dgg,\de)$-contact process are well-defined and the processes
started from a finite initial state are a.s.\ finite for all time.

For any $A\sub\La$, we let
\be\label{rhodef}
\rho(A):=\P[\eta^A_t\neq\emptyset\ \forall t\geq 0]
=\P[A\times\{0\}\leadsto\infty]
\ee
denote the survival probability of the $(\La,a,\de)$-contact process
started in $A$. Similarly, $\rho^\dgg$ denotes the survival
probability of the $(\La,a^\dgg,\de)$-contact process. Setting
\be\label{oveta}
\ov\eta_t:=\{i\in\La:-\infty\leadsto(i,t)\}\qquad(t\in\R)
\ee
defines a stationary $(\La,a,\de)$-contact process whose invariant law
\be\label{ovnu}
\ov\nu:=\P[\ov\eta_t\in\cdot\,]\qquad(t\in\R)
\ee
is uniquely characterized by
\be\label{nuchar}
\P\big[\ov\eta_0\cap A\neq\emptyset\big]=\rho^\dgg(A)
\qquad(A\in\Pc_{\rm fin}(\La)).
\ee
(To see this, note that the linear span of the functions
$B\mapsto1_{\{A\cap B=\emptyset\}}$ with $A\in\Pc_{\rm fin}(\La)$
forms an algebra that separates points, hence by the Stone-Weierstrass
theorem is dense in the space of continuous functions on $\Pc(\La)$.)
It is easy to see that $\ov\nu$ is {\em nontrivial}, i.e., gives zero
probability to the empty set, if and only if the
$(\La,a^\dgg,\de)$-contact process survives. Moreover, $\ov\nu$ is the limit
law of the process started with all sites infected, i.e., $\ov\nu$ is the
upper invariant law.

\subsection{Martingale problem}

We will need the fact that $(\La,a,\de)$-contact processes started in
finite initial states solve a martingale problem. Let
\be\label{Gdom}
\Si(\Pc_{\rm fin}(\La)):=\{f:\Pc_{\rm fin}(\La)\to\R:
|f(A)|\leq K|A|^k+M\mbox{ for some }K,M,k\geq 0\}.
\ee
denote the class of real functions on $\Pc_{\rm fin}(\La)$ of
polynomial growth. Given $\La,a$, and $\de$, define a linear operator
$G$ with domain $\Di(G):=\Si(\Pc_{\rm fin}(\La))$ by
\bc\label{Gdef}
Gf(A)&:=&\dis\sum_{ij}a(i,j)1_{\{i\in A\}}1_{\{j\not\in A\}}
\{f(A\cup\{j\})-f(A)\}\\[5pt]
&&\dis+\de\sum_i1_{\{i\in A\}}\{f(A\beh\{i\})-f(A)\}.
\ec
\bp{\bf(Martingale problem and moment estimate)}\label{P:mart}
The operator $G$ maps the space $\Si(\Pc_{\rm fin}(\La))$ into itself.
For each $f\in\Si(\Pc_{\rm fin}(\La))$ and $A\in\Pc_{\rm fin}(\La)$,
the process
\be\label{MP}
M_t:=f(\eta^A_t)-\int_0^tGf(\eta^A_s)\di s\qquad(t\geq 0)
\ee
is a martingale with respect to the filtration generated by $\eta^A$.
Moreover, setting $z^{\li k\re}:=\prod_{i=0}^{k-1}(z+i)$, one has
\be\label{momest}
\E\big[|\eta^A_t|^{\li k\re}\big]\leq|A|^{\li k\re}e^{k(|a|+(k-2)\de)t}
\qquad(A\in\Pc_{\rm fin}(\La),\ k\geq 1,\ t\geq 0).
\ee
\ep
{\bf Proof} Our proof follows the same lines as the proof of
\cite[Proposition~8]{AS05}. It is not hard to see that $G$ maps
$\Si(\Pc_{\rm fin}(\La))$ into itself.
\detail{If $f$ satisfies $|f(A)|\leq K|A|^k+M$ for some $K,M,k\geq 0$, then
\bc
|Gf(A)|&\leq&\dis\sum_{ij}a(i,j)1_{\{i\in A\}}1_{\{j\not\in A\}}
|f(A\cup\{j\})-f(A)|\\[5pt]
&&\dis+\de\sum_i1_{\{i\in A\}}|f(A\beh\{i\})-f(A)|\\[5pt]
&\leq&\dis\sum_{ij}a(i,j)1_{\{i\in A\}}
\big(K(|A|^{k+1}+|A|^k)+2M\big)\\[5pt]
&&\dis+\de\sum_i1_{\{i\in A\}}\big(K(|A|^k+|A|^{k-1})+2M\big)\\[5pt]
&\leq&\dis|a|\big(K(|A|^{k+1}+|A|^k)+2M\big)\\[5pt]
&&\dis+\de|A|\big(K(|A|^k+|A|^{k-1})+2M\big).
\ec}
Set $f_k(A):=|A|^{\li k\re}$. Then, using the fact that $z^{\li
  k\re}-(z-1)^{\li k\re}=kz^{\li k-1\re}$, we see that
\bc\label{Gmom}
\dis Gf_k(A)&=&\dis\sum_{ij}a(i,j)1_{\{i\in A\}}1_{\{j\not\in A\}}
\{(|A|+1)^{\li k\re}-|A|^{\li k\re}\}+\de\sum_i1_{\{i\in A\}}
\{(|A|-1)^{\li k\re}-|A|^{\li k\re}\},\\[5pt]
&\leq&\dis|a||A|\{(|A|+1)^{\li k\re}-|A|^{\li k\re}\}
-\de|A|\{|A|^{\li k\re}-(|A|-1)^{\li k\re}\}\\[5pt]
&=&\dis\big(|a|-\de\big)k|A|^{\li k\re}+\de k(k-1)|A|^{\li k-1\re}
\leq k\big(|a|+(k-2)\de\big)|A|^{\li k\re}=:K_k|A|^{\li k\re}.
\ec
\detail{Here we are using that $z^{\li
  k\re}-(z-1)^{\li k\re}=(z+k-1)z^{\li k-1\re}-(z-1)z^{\li k-1\re}=kz^{\li
  k-1\re}$, and $|a|z\{(z+1)^{\li k\re}-z^{\li k\re}\}-\de z\{z^{\li
  k\re}-(z-1)^{\li k\re}\}=|a|zk(z+1)^{\li k-1\re}-\de zkz^{\li k-1\re}=
  |a|kz^{\li k\re}-\de zk(z+1)^{\li k-1\re}
  +\de zk((z+1)^{\li k-1\re}-z^{\li k-1\re})
  =(|a|-\de)kz^{\li k\re}+\de zk(k-1)(z+1)^{\li k-2\re}
  =(|a|-\de)kz^{\li k\re}+\de k(k-1)z^{\li k-1\re}$.}
Define stopping times $\tau_N:=\inf\{t\geq 0:|\eta^A_t|\geq N\}$. The
stopped process $(\eta^A_{t\wedge\tau_N})_{t\geq 0}$ has bounded jump
rates, and therefore standard theory tells us that for each $N\geq 1$
and $f\in\Si(\Pc_{\rm fin}(\La))$, the process
\be\label{MNP}
M^N_t:=f(\eta^A_{t\wedge\tau_N})-\int_0^{t\wedge\tau_N}Gf(\eta^A_s)\di s
\qquad(t\geq 0)
\ee
is a martingale. Moreover, it easily follows from (\ref{Gmom}) that
\be\label{stopest}
\E\big[|\eta^A_{t\wedge\tau_N}|^{\li k\re}\big]\leq|A|^{\li k\re}e^{K_kt}
\qquad(k\geq 1,\ t\geq 0),
\ee
which in turn implies that $\P[|\eta^A_{t\wedge\tau_N}|\geq N]\to 0$ as
$N\to\infty$, and hence $\lim_{N\to\infty}\tau_N=\infty$. Using the fact that
$G$ maps $\Si(\Pc_{\rm fin}(\La))$ into itself and (\ref{stopest}) for some
sufficiently high $k$ (depending on $f$), one can show that for fixed $t\geq
0$, the random variables $(M^N_t)_{N\geq 1}$ are uniformly
integrable. Therefore, letting $N\to\infty$ in (\ref{MNP}), one finds that the
process in (\ref{MP}) is a martingale. Letting $N\to\infty$ in (\ref{stopest})
yields (\ref{momest}).\qed

\subsection{Covariance formula}\label{S:cov}

By Proposition~\ref{P:mart}, setting
\be\label{St}
S_tf(A):=\E[f(\eta^A_t)]
\qquad(f\in\Si(\Pc_{\rm fin}(\La)),\ A\in\Pc_{\rm fin}(\La))
\ee
defines a semigroup $(S_t)_{t\geq 0}$ of linear operators $S_t:\Si(\Pc_{\rm
  fin}(\La))\to\Si(\Pc_{\rm fin}(\La))$.  Let $\Mi$ be the class of
probability measures on $\Pc_{\rm fin}(\La)$ such that $\int|A|^k\mu(\di
A)<\infty$ for all $k\geq 1$. For $\mu\in\Mi$ and $f\in\Si(\Pc_{\rm
  fin}(\La))$, we write $\mu f:=\int f(A)\mu(\di A)$. Note that if
$(\eta_t)_{t\geq 0}$ is a $(\La,a,\de)$-contact process started in an initial
law $\P[\eta_0\in\cdot\,]=:\mu\in\Mi$, then $\P[\eta_t\in\cdot\,]\in\Mi$ for
all $t\geq 0$ and $\int\P[\eta_t\in\di A]f(A)=\mu S_tf$. For this reason, we
use the notation $\mu S_t:=\P[\eta_t\in\cdot\,]$ $(t\geq 0)$ to denote the law
of $\eta_t$. For any $\mu\in\Mi$ and $f,g\in\Si(\Pc_{\rm fin}(\La))$, we let
\be
\cov_\mu(f,g):=\mu(fg)-(\mu f)(\mu g)
\ee
denote the covariance of $f$ and $g$ under $\mu$, which is always
finite.
\bp{\bf(Covariance formula)}\label{P:cov}
For $f,g\in\Si(\Pc_{\rm fin}(\La))$, let
\be\label{Gadef}
\Ga(f,g):=\ffrac{1}{2}\big[G(fg)-(Gf)g-f(Gg)\big]
\qquad(f,g\in\Si(\Pc_{\rm fin}(\La))).
\ee
Then, for any $\mu\in\Mi$ and $f,g\in\Si(\Pc_{\rm fin}(\La))$, one has
\be
\cov_{\mu S_t}(f,g)=\cov_\mu(S_tf,S_tg)
+2\!\int_0^t\mu S_{t-s}\Ga(S_sf,S_sg)\,\di s\qquad(t\geq 0).
\ee
\ep
{\bf Proof} Set
\be
H(s,t,u):=S_s\big((S_tf)(S_ug)\big).
\ee
We claim that
\bc\label{Hstu}
\dis\dif{s}H(s,t,u)&=&\dis S_sG\big((S_tf)(S_ug)\big),\\[5pt]
\dis\dif{t}H(s,t,u)&=&\dis S_s\big((GS_tf)(S_ug)\big),\\[5pt]
\dis\dif{u}H(s,t,u)&=&\dis S_s\big((S_tf)(GS_ug)\big).
\ec
It follows that
\be
\dif{t}H(t,T-t,T-t)=2S_t\Ga(S_{T-t}f,S_{T-t}g),
\ee
and therefore
\be\ba{l}
\dis\cov_{\mu S_T}(f,g)-\cov_\mu(S_Tf,S_Tg)\\[5pt]
\dis\quad=\big(\mu S_T(fg)-(\mu S_Tf)(\mu S_Tg)\big)
-\big(\mu((S_Tf)(S_Tg))-(\mu S_Tf)(\mu S_Tg)\big)\\[5pt]
\dis\quad=\mu\big(S_T(fg)-(S_Tf)(S_Tg)\big)\\[-2pt]
\dis\quad=\mu\big(H(T,0,0)-H(0,T,T)\big)
=2\int_0^t\mu S_t\Ga(S_{T-t}f,S_{T-t}g)\,\di t.
\ec
These calculations are standard. However, in order to verify
(\ref{Hstu}), we must use some special properties of our model.
Let us say that a sequence of functions $f_n\in\Si(\Pc_{\rm
fin}(\La))$ converges `nicely' to a limit $f$, if $f_n\to f$ pointwise
and there exists $K,M,k\geq 0$ such that $|f_n(A)|\leq K|A|^k+M$ for
all $n$. By (\ref{momest}) and dominated convergence, if $f_n\to f$
`nicely', then $S_tf_n\to S_t f$ `nicely', for each $t\geq 0$.  Note
also that if $f_n,f,g\in\Si(\Pc_{\rm fin}(\La))$ and $f_n\to f$
`nicely', then $f_ng\to fg$ `nicely'. We claim that for each
$f\in\Si(\Pc_{\rm fin}(\La))$,
\be\label{GSt}
\lim_{t\to 0}t^{-1}(S_tf-f)=Gf,
\ee
where the convergence happens `nicely'. Indeed, by Proposition~\ref{P:mart},
\be
t^{-1}\big(S_tf(A)-f(A)\big)=t^{-1}\!\int_0^t\E\big[(Gf)(\eta^A_s)\big]\di s
\astol{t\to 0}Gf(A),
\ee
where the `niceness' of the convergence follows from (\ref{momest}) and the
fact that $Gf\in\Si(\Pc_{\rm fin}(\La))$. It follows from (\ref{GSt}) that for
each $f\in\Si(\Pc_{\rm fin}(\La))$ and $t\geq 0$,
\bc\label{difSt}
\dis\dif{t}S_tf
&=&\dis\lim_{\eps\to0}\eps^{-1}(S_\eps-1)S_tf=GS_tf\\[5pt]
&=&\dis\lim_{\eps\to0}\eps^{-1}S_t(S_\eps-1)f=S_tGf,
\ec
where $1$ denotes the identity operator. Using (\ref{difSt}) and the
properties of `nice' convergence, (\ref{Hstu}) follows readily.\qed

\section{The exponential growth rate}\label{S:exp}

\subsection{Basic facts}\label{S:basic}

{\bf Proof of Lemma~\ref{L:exprate}} Let us write
\be\label{pitdef}
\expi_t(A):=\E\big[|\eta^A_t|\big]\qquad(A\in\Pc_{\rm fin}(\La),\ t\geq 0).
\ee
We start by showing that
\be\label{pimult}
\expi_{s+t}(\{0\})\leq\expi_s(\{0\})\expi_t(\{0\})\qquad(s,t\geq 0).
\ee
If $\eta^A$ and $\eta^B$ are defined usng the same graphical
representation, then $\eta^A_t\cup\eta^B_t=\eta^{A\cup
B}_t$. Therefore,
\be\label{pixpi}
\E\big[|\eta^A_t|\big]=\E\Big[\big|\bigcup_{i\in A}\eta^{\{i\}}_t\big|\Big]
\leq\sum_{i\in A}\E\big[|\eta^{\{i\}}_t|\big]=|A|\E\big[|\eta^{\{0\}}_t|\big],
\ee
where in the last step we have used shift invariance. As a consequence,
\be\label{submult}
\expi_{s+t}(\{0\})=\int\P[\eta^{\{0\}}_s\in\di A]\E\big[|\eta^A_t|\big]
\leq\int\P[\eta^{\{0\}}_s\in\di A]|A|\E\big[|\eta^{\{0\}}_t|\big]
=\expi_s(\{0\})\expi_t(\{0\}).
\ee
This proves (\ref{pimult}). It follows that $t\mapsto\log\expi_t(\{0\})$
is subadditive and therefore, by \cite[Theorem~B.22]{Lig99}, the limit
\be\label{rdef}
\lim_{t\to\infty}\,\ffrac{1}{t}\log\expi_t(\{0\})
=\inf_{t>0}\ffrac{1}{t}\log\expi_t(\{0\})=:r\in[-\infty,\infty)
\ee
exists. By monotonicity and (\ref{pixpi}),
\be\label{piA}
\expi_t(\{0\})\leq\expi_t(A)\leq|A|\expi_t(\{0\})
\qquad(A\in\Pc_{\rm fin}(\La)).
\ee
Taking logarithms, dividing by $t$, and letting $t\to\infty$ we arrive
at (\ref{expr}). Since $\eta$ can be bounded from below by a simple
death process and from above by a branching process, one has
\be
e^{-\de t}\leq\E\big[|\eta^{\{0\}}_t|\big]\leq e^{(|a|-\de)t}\qquad(t\geq 0),
\ee
which implies that $-\de\leq r\leq|a|-\de$.\qed

\noi
{\bf Proof of Theorem~\ref{T:main}~(a)} By duality (formula (\ref{dual}))
and shift invariance,
\bc\label{Esym}
\dis\E\big[|\eta^{\{0\}}_t|\big]
&=&\dis\sum_i\P\big[\eta^{\{0\}}_t\cap\{i\}\neq\emptyset\big]
=\sum_i\P\big[\{0\}\cap\eta^{\dgg\,\{i\}}_t\neq\emptyset\big]\\[5pt]
&=&\dis\sum_i\P\big[\{i^{-1}\}\cap\eta^{\dgg\,\{0\}}_t\neq\emptyset\big]
=\E\big[|\eta^{\dgg\,\{0\}}_t|\big],
\ec
which implies that $r(\La,a,\de)=r(\La,a^\dgg,\de)$.\qed

\noi
{\bf Proof of Theorem~\ref{T:main}~(b)} Fix a countable group $\La$ and
infection rates $a$ satisfying (\ref{assum}), and for each $\de\geq 0$, write
$\pi(\de,t):=\E\big[|\eta^{\{0\}}_t|\big]$, where $\eta^{\{0\}}_t$ is the
$(\La,a,\de)$-contact process. For $0\leq\de<\ti\de$, consider the graphical
representations (see Section~\ref{S:graph}) of the $(\La,a,\de)$- and
$(\La,a,\ti\de)$-contact processes, defined by Poisson processes $(\oo^{\rm
  r},\oo^{\rm i})$ and $(\ti\oo^{\rm r},\ti\oo^{\rm i})$, respectively. We may
couple these graphical representations such that $\oo^{\rm i}=\ti\oo^{\rm i}$
and $\oo^{\rm r}\sub\ti\oo^{\rm r}$, where $\ti\oo^{\rm r}\beh\oo^{\rm r}$ is
an Poisson point process with intensity $\ti\de-\de$, independent of $\oo^{\rm
  i}$ and $\oo^{\rm r}$. Write $\leadsto$ and $\ti\leadsto$ to indicate the
existence of an open path in the graphical representations for $\de$ and
$\ti\de$, respectively. Then, if $\ti\de-\de$ is small, then for each $t\geq
0$,
\be\ba{l}\label{lipi}
\dis\pi(\ti\de,t)=\sum_i\P[(0,0)\ti\leadsto(i,t)]\\[5pt]
\dis\quad=\sum_i\P[(0,0)\leadsto(i,t)]\\[5pt]
\dis\quad\phantom{=}-(\ti\de-\de)\int_0^t\sum_{ij}
\P\big[(0,0)\leadsto(j,s)\leadsto(i,t)\mbox{ and there exists}\\[-10pt]
\dis\quad\phantom{=+(\ti\de-\de)\int_0^t\sum_{ij}\P\big[}
\mbox{no $k\neq j$ such that }(0,0)\leadsto(k,s)\leadsto(i,t)\big]\di s
+O((\ti\de-\de)^2),
\ec
where the terms order $(\ti\de-\de)^2$ come from events where two or
more recovery symbols in $\ti\oo^{\rm r}\beh\oo^{\rm r}$ are needed to
block all paths from $(0,0)$ to $(i,t)$. Dividing by $\ti\de-\de$ and
letting $\ti\de\to\de$ yields
\be\ba{l}
\dis\dif{\de}\dis\pi(\ti\de,t)\\[5pt]
\dis\quad=-\int_0^t\sum_i\P\big[\exists j\mbox{ s.t.\ }
(0,0)\leadsto(j,s)\leadsto(i,t)\mbox{ and }
\not\exists k\neq j\mbox{ s.t.\ }(0,0)\leadsto(k,s)\leadsto(i,t)\big]\di s,
\ec
which is an analogue of what is known as Russo's formula in percolation.
Since
\be\ba{l}
\dis\sum_i\P\big[\exists j\mbox{ s.t.\ }
(0,0)\leadsto(j,s)\leadsto(i,t)\mbox{ and }
\not\exists k\neq j\mbox{ s.t.\ }(0,0)\leadsto(k,s)\leadsto(i,t)\big]\\[5pt]
\dis\quad\leq\sum_i\P[(0,0)\leadsto(i,t)]=\pi(\de,t),
\ec
it follows that $0\leq-\dif{\de}\pi(\de,t)\leq t\pi(\de,t)$ $(t\geq 0)$,
and therefore
\be
0\leq-\dif{\de}\ffrac{1}{t}\log\pi(\de,t)\leq 1.
\ee
Taking the limit $t\to\infty$, using (\ref{rdef}), the claims follow.
Note that by Lemma~\ref{L:exprate}, $-\de\leq r\leq|a|-\de$, so letting
$\de\to\infty$ we see that the Lipschitz constant $1$ is optimal.\qed

\noi
{\bf Proof of Theorem~\ref{T:main}~(c)} If the $(\La,a,\de)$-contact
process survives, then
\be
\E\big[|\eta^{\{0\}}_t|\big]\geq\P[\eta^{\{0\}}_t\neq\emptyset]
\asto{t}\P[\eta^{\{0\}}_s\neq\emptyset\ \forall s\geq 0]>0,
\ee
which implies that $r\geq 0$.\qed

\subsection{Eigenmeasures}

Recall that a measure $\mu$ on a locally compact space is called
locally finite if $\mu(K)<\infty$ for all compact sets $K$. We need a
few basic facts about locally finite measures on $\Pc_+(\La)$.
\bl{\bf(Locally finite measures)}\label{L:locfin}
Let $\mu$ be a measure on $\Pc_+(\La)$. Then the following statements are
equivalent: 1.\ $\mu$ is locally finite. 2.\ $\int\mu(\di A)1_{\{i\in
A\}}<\infty$ for all $i\in\La$. 3.\ $\int\mu(\di A)1_{\{A\cap
B\neq\emptyset\}}<\infty$ for all $B\in\Pc_{\rm fin}(\La)$.
\el
{\bf Proof} We will prove that 1$\volgt$3$\volgt$2$\volgt$1. For each
$B\in\Pc_{\rm fin}(\La)$, the set $\Qi(B):=\{A\sub\La:A\cap B\neq\emptyset\}$
is a compact subset of $\Pc_+(\La)$, and $A\mapsto1_{\{A\cap
  B\neq\emptyset\}}$ is a continuous function with compact support
$\Qi(B)$. It follows that any locally finite measure $\mu$ satisfies
$\int\mu(\di A)1_{\{A\cap B\neq\emptyset\}}<\infty$ for each $B\in\Pc_{\rm
  fin}(\La)$. In particular, setting $B=\{i\}$ this implies that $\int\mu(\di
A)1_{\{i\in A\}}<\infty$ for each $i\in\La$. This proves the implications
1$\volgt$3$\volgt$2. To see that 2 implies 1, let $\De_n\sub\La$ be finite
sets increasing to $\La$. Then the $\Qi(\De_n)$ increase to $\Pc_+(\La)$ and,
since the $\Qi(\De_n)$ are open sets, each compact subset of
$\Pc_+(\La)$ is contained in some $\Qi(\De_n)$. Therefore, since
$\mu(\Qi(\De_n))=\int\mu(\di A)1_{\{A\cap\De_n\neq\emptyset\}}
\leq\sum_{i\in\De_n}\int\mu(\di A)1_{\{i\in A\}}<\infty$ for each $n$,
the measure $\mu$ is locally finite.\qed

\noi
We equip the space of locally finite measures on $\Pc_+(\La)$ with the
vague topology, i.e., we say that a sequence of locally finite
measures $\mu_n$ on $\Pc_+(\La)$ converges vaguely to a limit $\mu$,
denoted as $\mu_n\Rightarrow\mu$, if $\int\mu_n(\di
A)f(A)\to\int\mu(\di A)f(A)$ for each continuous compactly supported
real function $f$ on $\Pc_+(\La)$.
\bl{\bf(Vague convergence)}\label{L:vague}
Let $\mu_n,\mu$ be locally finite measures on $\Pc_+(\La)$. Then the
$\mu_n$ converge vaguely to $\mu$ if and only if $\int\mu_n(\di
A)1_{\{A\cap B\neq\emptyset\}}\to\int\mu(\di A)1_{\{A\cap
B\neq\emptyset\}}$ for each $B\in\Pc_{\rm fin}(\La)$. The sequence
$\mu_n$ is relatively compact in the topology of vague convergence if
and only if $\sup_n\int\mu_n(\di A)1_{\{A\cap
B\neq\emptyset\}}<\infty$ for each $B\in\Pc_{\rm fin}(\La)$.
\el
{\bf Proof} Since for each $B\in\Pc_{\rm fin}(\La)$, the function
$A\mapsto1_{\{A\cap B\neq\emptyset\}}$ is continuous with compact
support, the conditions for convergence and relative compactness given
above are clearly necessary. To see that they are also sufficient, let
$\De_m\sub\La$ be finite sets increasing to $\La$ and set
$f_m(A):=1_{\{A\cap\De_n\neq\emptyset\}}$. Then the $f_m$ are
continuous, nonnegative functions with compact supports increasing to
$\Pc_+(\La)$. It follows that $\mu_n$ converges vaguely to $\mu$ if
and only if for each $m$ the weighted measures $f_m(A)\mu_n(\di A)$
converge weakly to $f_m(A)\mu(\di A)$. Now if $\sup_n\int\mu_n(\di
A)1_{\{A\cap B\neq\emptyset\}}<\infty$ for each $B\in\Pc_{\rm
fin}(\La)$, then by a diagonal argument each subsequence of the
$\mu_n$ contains a further subsequence such that $f_m(A)\mu_n(\di A)$
converges weakly for each $m$, hence the $\mu_n$ converge vaguely. The
linear span of the functions $B\mapsto1_{\{A\cap B=\emptyset\}}$ with
$A\in\Pc_{\rm fin}(\La)$ forms an algebra that separates points, hence
by the Stone-Weierstrass theorem is dense in the space of continuous
functions on $\Pc(\La)$. It follows that $\mu_n$ converges vaguely to
$\mu$ if and only if $\int f_m(A)\mu_n(\di A)1_{\{A\cap
B\neq\emptyset\}}$ converges to $\int f_m(A)\mu(\di A)1_{\{A\cap
B\neq\emptyset\}}$ for each $m$ and for each $B\in\Pc_{\rm
fin}(\La)$. Since $f_m(A)1_{\{A\cap
B\neq\emptyset\}}=1_{\{A\cap\De_n\neq\emptyset\}} +1_{\{A\cap
B\neq\emptyset\}}-1_{\{A\cap(B\cup\De_n)\neq\emptyset\}}$, this is in
turn implied by the condition that $\int\mu_n(\di A)1_{\{A\cap
B\neq\emptyset\}}\to\int\mu(\di A)1_{\{A\cap B\neq\emptyset\}}$ for
each $B\in\Pc_{\rm fin}(\La)$.\qed

\detail{$f_m(A)1_{\{A\cap B\neq\emptyset\}}
=(1-1_{\{A\cap\De_n=\emptyset\}})(1-1_{\{A\cap B=\emptyset\}})
=1-1_{\{A\cap\De_n=\emptyset\}}-1_{\{A\cap B=\emptyset\}}
+1_{\{A\cap(B\cup\De_n)=\emptyset\}}
=1_{\{A\cap\De_n\neq\emptyset\}}+1_{\{A\cap B\neq\emptyset\}}
-1_{\{A\cap(B\cup\De_n)\neq\emptyset\}}$.}

\noi
The next lemma guarantees that expressions as in the left-hand
side of (\ref{eigen}) are well-defined and yield a homogeneous,
locally finite measure on $\Pc_+(\La)$
\bl{\bf(Evolution of locally finite measures)}\label{L:infev}
If $\mu$ is a homogeneous, locally finite measure on $\Pc_+(\La)$,
then for each $t\geq 0$, the measure $\int\mu(\di
A)\P[\eta^A_t\in\cdot\,]|_{\Pc_+(\La)}$ is homogeneous and locally finite on
$\Pc_+(\La)$. If $\mu_n$ are homogeneous, locally finite measures on
$\Pc_+(\La)$ converging vaguely to a limit $\mu$, then
\be\label{mueta}
\int\mu_n(\di A)\P[\eta^A_t\in\cdot\,]\big|_{\Pc_+(\La)}
\Asto{n}\int\mu(\di A)\P[\eta^A_t\in\cdot\,]\big|_{\Pc_+(\La)}\qquad(t\geq 0).
\ee
\el
{\bf Proof} We start by observing that any homogeneous, locally finite
measure $\mu$ on $\Pc_+(\La)$ satisfies
\be\label{homest}
\int\mu(\di A)1_{\{A\cap B\neq\emptyset\}}
\leq\sum_{i\in B}\int\mu(\di A)1_{\{i\in A\}}
=|B|\int\mu(\di A)1_{\{0\in A\}}<\infty.
\ee
Using duality (see (\ref{dual})), it follows that
\be\ba{l}\label{ABest}
\dis\int\mu(\di A)\int\P[\eta^A_t\in\di C]1_{\{C\cap B\neq\emptyset\}}
=\int\mu(\di A)\P[\eta^A_t\cap B\neq\emptyset]
=\int\mu(\di A)\P[A\cap\eta^{\dgg\,B}_t\neq\emptyset]\\[5pt]
\dis\quad=\int\P[\eta^{\dgg\,B}_t\in\di C]
\int\mu(\di A)1_{\{A\cap C\neq\emptyset\}}
\leq\int\P[\eta^{\dgg\,B}_t\in\di C]\:|C|\int\mu(\di A)1_{\{0\in A\}}\\[5pt]
\dis\quad=\E\big[|\eta^{\dgg\,B}_t|\big]\int\mu(\di A)1_{\{0\in A\}}<\infty
\qquad\qquad(B\in\Pc_{\rm fin}(\La)).
\ec
By Lemma~\ref{L:locfin}, it follows that the measure $\int\mu_n(\di
A)\P[\eta^A_t\in\cdot\,]|_{\Pc_+(\La)}$ is locally finite. It is obviously
homogeneous. Now if $\mu_n$ are homogeneous, locally finite measures
on $\Pc_+(\La)$ converging vaguely to a limit $\mu$, then, by the first three
equalities in (\ref{ABest}),
\be
\dis\int\mu_n(\di A)\int\P[\eta^A_t\in\di C]1_{\{C\cap B\neq\emptyset\}}
=\int\P[\eta^{\dgg\,B}_t\in\di C]\:\int\mu_n(\di A) 1_{\{A\cap
C\neq\emptyset\}},
\ee
for each $B\in\Pc_{\rm fin}(\La)$, and this quantity converges to the analogue
quantity for $\mu$ by dominated convergence, using (\ref{fiex}), the
estimate (\ref{homest}), and the fact that the $\int\mu_n(\di A)1_{\{0\in A\}}$
are uniformly bounded in $n$ since they converge. Applying
Lemma~\ref{L:vague}, we arrive at (\ref{mueta}).\qed

\noi
{\bf Proof of Proposition~\ref{P:eigen}} It suffices to prove, for
each $(\La,a,\de)$-contact process, the following three claims:
\begin{enumerate}
\item There exists a homogeneous eigenmeasure of the
$(\La,a,\de)$-contact process, with eigenvalue $r=r(\La,a,\de)$.

\item If $\la$ is the eigenvalue of a homogeneous eigenmeasure
  of the $(\La,a,\de)$-contact process, then $\la\leq r$.

\item The set $\Ei(\La,a,\de)$ is closed and bounded from below.
\end{enumerate}
We start with claim~1. Define (by Lemma~\ref{L:infev} applied to
$\mu=\sum_i\de_{\{i\}}$) homogeneous, locally finite measures $\mu_t$ on
$\Pc_+(\La)$ by
\be\label{mut}
\mu_t:=\sum_i\P[\eta^{\{i\}}_t\in\cdot\,]\big|_{\Pc_+(\La)}\qquad(t\geq 0).
\ee
Let $\hat\mu_\la$ be the Laplace transform of $(\mu_t)_{t\geq 0}$, i.e.,
\be\label{Lap}
\hat\mu_\la:=\int_0^\infty\mu_t\:e^{-\la t}\di t\qquad(\la>r).
\ee
We claim that the measures $\hat\mu_\la$ are locally finite and, properly
renormalized, relatively compact in the topology of vague convergence, and
that each subsequential limit as $\la\down r$ is a homogeneous
eigenmeasure of the $(\La,a,\de)$-contact process, with eigenvalue~$r$.

Note that by duality (see (\ref{dual})),
\be\label{mudu}
\int\mu_t(\di A)1_{\{A\cap B\neq\emptyset\}}
=\sum_i\P[\eta^{\{i\}}_t\cap B\neq\emptyset]
=\sum_i\P[i\in\eta^{\dgg\,B}_t]
=\E\big[|\eta^{\dgg\,B}_t|\big]=\pi^\dgg_t(B)
\ee
$(t\geq 0,\ B\in\Pc_{\rm fin}(\La))$, where $\pi^\dgg_t(A)$ is defined
in analogy with (\ref{pitdef}) for the $(\La,a^\dgg,\de)$-contact
process. It follows that
\be\label{hatmudu}
\int\hat\mu_\la(\di A)1_{\{A\cap B\neq\emptyset\}}=\hat\pi^\dgg_\la(B),
\ee
where
\be
\hat\pi^\dgg_\la(A):=
\int_0^\infty\pi^\dgg_t(A)\:e^{-\la t}\di t
\qquad(\la>r,\ A\in\Pc_{\rm fin}(\La)).
\ee
By (\ref{rdef}) and Theorem~\ref{T:main}~(a), which has been proved in
Section~\ref{S:basic}, for every $\eps>0$, there exists a
$T_\eps<\infty$ such that
\be\label{twosid}
e^{rt}\leq\pi^\dgg_t(\{0\})\leq e^{(r+\eps)t}\qquad(t\geq T_\eps).
\ee
It follows from the upper bound in (\ref{twosid}) that
$\hat\pi^\dgg_\la(\{0\})<\infty$ for all $\la>r$. Hence, by
(\ref{hatmudu}) and Lemma~\ref{L:locfin}, the measures $\hat\mu_\la$
are locally finite for each $\la>r$. The lower bound in (\ref{twosid})
and monotone convergence moreover show that
\be\label{piinf}
\lim_{\la\down r}\hat\pi^\dgg_\la(\{0\})
=\lim_{\la\down r}\int_0^\infty\pi^\dgg_t(\{0\})\:e^{-\la t}\di t
=\int_0^\infty\pi^\dgg_t(\{0\})\:e^{-rt}\di t=\infty.
\ee
Set $\ov\mu_\la:=\hat\pi^\dgg_\la(\{0\})^{-1}\hat\mu_\la$. Then for each
$\la>r$, the measure $\ov\mu_\la$ is homogenous, locally finite, and
normalized such that $\int\ov\mu_\la(\di A)1_{\{0\in
A\}}=1$. Therefore, by Lemma~\ref{L:vague} and the estimate
(\ref{homest}), the measures $\ov\mu_\la$ are relatively compact in
the topology of vague convergence as $\la\down r$. Choose $\la_n\down
r$ such that $\ov\mu_{\la_n}\Rightarrow\ov\mu_r$ for some homogenous,
locally finite $\ov\mu_r$. Obviously $\int\ov\mu_r(\di A)1_{\{0\in
A\}}=1$ so $\ov\mu_r$ is nonzero. Then, filling in our definitions, using
Lemma~\ref{L:infev} and the Markov property of the contact process,
\be\ba{l}
\dis\int\ov\mu_r(\di A)\P[\eta^A_t\in\cdot\,]\big|_{\Pc_+(\La)}
=\lim_{n\to\infty}\hat\pi^\dgg_{\la_n}(\{0\})^{-1}
\int_0^\infty\!\!e^{-\la_n s}\di s\:\sum_i\int\P[\eta^{\{i\}}_s\in\di A]
\P[\eta^A_t\in\cdot\,]\big|_{\Pc_+(\La)}\\[5pt]
\dis\qquad=e^{rt}\lim_{n\to\infty}\hat\pi^\dgg_{\la_n}(\{0\})^{-1}
\int_0^\infty\!\!e^{-\la_n(s+t)}\di s\:\sum_i
\P[\eta^{\{i\}}_{s+t}\in\cdot\,]\big|_{\Pc_+(\La)}\\[5pt]
\dis\qquad=e^{rt}\Big(\ov\mu_r
-\lim_{n\to\infty}\hat\pi^\dgg_{\la_n}(\{0\})^{-1}
\int_0^t\!\!e^{-\la_n s}\di s\:\sum_i
\P[\eta^{\{i\}}_s\in\cdot\,]\big|_{\Pc_+(\La)}\Big)
=e^{rt}\ov\mu_r,
\ec
where in the last step we have used (\ref{piinf}). This shows that
$\ov\mu_r$ is an eigenmeasure with eigenvalue $r$.

To prove claim~2, we observe that if $\mu$ is a homogeneous
eigenmeasure with eigenvalue $\la$, then by duality (see (\ref{dual})) and
(\ref{homest}),
\be\ba{l}
\dis e^{\la t}\int\mu(\di A)1_{\{0\in A\}}
=\int\mu(\di A)\P[0\in\eta^A_t]\\[5pt]
\dis=\int\mu(\di A)\P[\eta^{\dgg\,\{0\}}_t\cap A\neq\emptyset]
\leq\E\big[|\eta^{\dgg\,\{0\}}_t|\big]\int\mu(\di A)1_{\{0\in A\}}.
\ec
By (\ref{twosid}), for each $\eps>0$ we can choose $t$ large enough
such that $\E\big[|\eta^{\dgg\,\{0\}}_t|\big]\leq
e^{(r+\eps)t}$. Since $\int\mu(\di A)1_{\{0\in A\}}>0$, we may divide
by it, yielding $e^{\la t}\leq e^{(r+\eps)t}$, which implies $\la\leq
r+\eps$. Since $\eps>0$ is arbitrary, it follows that $\la\leq r$.

To prove claim~3, finally, we observe that since we may estimate a
contact process from below by a simple death process, for any
homogeneous eigenmeasure $\mu$ with eigenvalue $\la$, one has
\be
e^{\la t}\int\mu(\di A)1_{\{0\in A\}}
=\int\mu(\di A)\P[0\in\eta^A_t]\geq e^{-\de t}\int\mu(\di A)1_{\{0\in A\}},
\ee
which shows that $\Ei(\La,a,\de)\sub[-\de,\infty)$. To show that
$\Ei(\La,a,\de)$ is closed, assume that $\la_n\in\Ei(\La,a,\de)$ and
$\la_n\to\la$. Then we can find homogeneous
eigenmeasures $\mu_n$ with eigenvalues $\la_n$. Normalizing such that
$\int\mu_n(\di A)1_{\{0\in A\}}=1$, using Lemma~\ref{L:vague} and
(\ref{homest}), we see that the sequence $\mu_n$ is relatively compact
in the topology of vague convergence, hence has a subsequential limit
$\mu$, which by Lemma~\ref{L:infev} is a homogeneous
eigenmeasure with eigenvalue $\la$.\qed

\noi
The proof of Proposition~\ref{P:eigen} yields a useful corollary.
\bcor{\bf(Convergence to eigenmeasure)}\label{C:conv}
Let $\mu_t$ be defined as in (\ref{mut}). Then the measures
\be
\hat\pi^\dgg_\la(\{0\})^{-1}\int_0^\infty\mu_t\:e^{-\la t}\di t
\ee
are relatively compact as $\la\down r$ in the topology of vague
convergence, and each subsequential limit as $\la\down r$ is a
homogeneous eigenmeasure of the $(\La,a,\de)$-contact
process, with eigenvalue $r(\La,a,\de)$.
\ecor
{\bf Remark} It seems intuitively plausible that the measures $\mu_t$,
suitably rescaled, converge as $t\to\infty$ to a vague limit, which by
Corollary~\ref{C:conv} then has to be an eigenmeasure with eigenvalue $r$.
Indeed, it seems plausible that this convergence is monotone, in a suitable
sense, and that these eigenmeasures are the `lowest' possible eigenmeasures,
in a suitable stochastic order. Should these conjectures be correct, then
these eigenmeasures are quite similar to the `second lowest invariant measure'
{f}rom \cite{SS97,SS99}, by which they are inspired. I do not know if these
conjectures are correct, or even what kind of stochastic order one should choose
here; I spent quite a bit of time in vain trying to prove that the measures
$\mu_t$ conditioned on the event that the origin is infected, are
stochastically increasing in time (in the usual stochastic order).

\subsection{Proof of Theorem~\ref{T:zero}}

We start with a preparatory lemma. We say that a function $f:\Pc_{\rm
fin}(\La)\to\R$ is {\em shift-invariant} if $f(iA)=f(A)$ for all
$i\in\La$, {\em monotone} if $A\sub B$ implies $f(A)\leq f(B)$, and
{\em subadditive} if $f(A\cup B)\leq f(A)+f(B)$, for all
$A,B\in\Pc_{\rm fin}(\La)$. Recall the definition of the generator $G$
of the $(\La,a,\de)$-contact process from (\ref{Gdef}). We define
$G^\dgg$ analogously, for the $(\La,a^\dgg,\de)$-contact process.
\bl{\bf(Eigenmeasures and harmonic functions)}\label{L:eighar}
If $\mu$ is a homogeneous eigenmeasure of the
$(\La,a,\de)$-contact process with eigenvalue $\la$, then
\be
v(A):=\int\mu(\di B)1_{\{A\cap B\neq\emptyset\}}\qquad(A\in\Pc_{\rm fin}(\La))
\ee
is a shift-invariant, monotone, subadditive function such that
$v(\emptyset)=0$, $v(A)>0$ for any $\emptyset\neq A\in\Pc_{\rm fin}(\La)$,
$v\in\Si(\Pc_{\rm fin}(\La))$, and $G^\dgg v=\la v$.
\el
{\bf Proof} The function $v$ is obviously shift-invariant, monotone, and
satisfies $v(\emptyset)=0$. Since $\mu$ is homogeneous and nonzero, $v(A)>0$
for any $\emptyset\neq A\in\Pc_{\rm fin}(\La)$. The function $v$ is
subadditive since $1_{\{(A\cup A')\cap B\neq\emptyset\}}\leq1_{\{A\cap
  B\neq\emptyset\}}+1_{\{A'\cap B\neq\emptyset\}}$. Subadditivity and
shift-invariance imply that $v(A)\leq v(\{0\})|A|$, so certainly
$v\in\Si(\Pc_{\rm fin}(\La))$.

To see that $G^\dgg v=\la v$, observe that by duality
(see (\ref{dual}))
\be\ba{l}\label{semigr}
\dis\E[v(\eta^{\dgg\,A}_t)]=\int\mu(\di B)\P[\eta^{\dgg\,A}_t\cap B\neq\emptyset]
=\int\mu(\di B)\P[A\cap\eta^B_t\neq\emptyset]\\[5pt]
\dis\quad=\int\mu(\di B)\P[\eta^B_t\in\di C]1_{\{A\cap C\neq\emptyset\}}
=e^{\la t}\int\mu(\di C)1_{\{A\cap C\neq\emptyset\}}=e^{\la t}v(A).
\ec
Let $(S^\dgg_t)_{t\geq 0}$ denote the semigroup of the
$(\La,a^\dgg,\de)$-contact process. Recall from Section~\ref{S:cov} that
$S^\dgg_t$ maps the space $\Si(\Pc_{\rm fin}(\La))$ into itself. Then
(\ref{semigr}) says that $S^\dgg_tv=e^{\la t}v$, and hence, by (\ref{GSt}),
$G^\dgg v=\lim_{\eps\to 0}\eps^{-1}(S^\dgg_\eps v-v)=\la v$.\qed

\noi
Recall the definition of the survival probability $\rho(A)$ from
(\ref{rhodef}). Theorem~\ref{T:zero} follows from the following,
stronger result.
\bp{\bf(Shift invariant monotone harmonic functions)}\label{P:monhar}
Assume that the infection rates satisfy the irreducibility condition
(\ref{irr}) and that the $(\La,a,\de)$-contact process survives. Assume that
$f:\Pc_{\rm fin}(\La)\to\R$ is shift invariant, monotone, $f(\emptyset)=0$,
$f\in\Si(\Pc_{\rm fin}(\La))$, and $Gf=0$. Then there exists a constant
$c\geq 0$ such that $f=c\rho$.
\ep
Before we prove this, we first show how this implies Theorem~\ref{T:zero}.\med

\noi
{\bf Proof of Theorem~\ref{T:zero}} Let $\mu$ be a homogeneous
eigenmeasure of the $(\La,a,\de)$-contact process with eigenvalue zero, and
let $v(A):=\int\mu(\di B)1_{\{A\cap B\neq\emptyset\}}$. By
Lemma~\ref{L:eighar}, $v$ is shift invariant, monotone, $v(\emptyset)=0$,
$v\in\Si(\Pc_{\rm fin}(\La))$, and $G^\dgg v=0$. By assumption, the upper
invariant measure $\ov\nu$ of the $(\La,a,\de)$-contact process is nontrivial,
hence the $(\La,a^\dgg,\de)$-contact process survives, so by
Proposition~\ref{P:monhar}, $v=c\rho^\dgg$ for some $c\geq 0$, where
$\rho^\dgg$ denotes the survival probability of the $(\La,a^\dgg,\de)$-contact
process. By the characterization of the upper invariant measure in
(\ref{nuchar}), it follows that $\mu=c\ov\nu$.\qed

\noi
In order to prove Proposition~\ref{P:monhar}, we need one more lemma.
\bl{\bf(Eventual domination of finite configurations)}\label{L:evdom}
Assume that the infection rates satisfy the irreducibility condition
(\ref{irr}) and that the $(\La,a,\de)$-contact process survives. Then
\be\label{evdom}
\lim_{t\to\infty}\P\big[\exists i\in\La\mbox{ s.t.\ }\eta^A_t\supset iB\,\big|
\,\eta^A_t\neq\emptyset]=1
\qquad(A,B\in\Pc_{\rm fin}(\La),\ A\neq\emptyset).
\ee
\el
Formula (\ref{evdom}) says that $\eta$ exhibits a form of `extinction
versus unbounded growth'. More precisely, either $\eta_t$ gets extinct
or $\eta_t$ is eventually larger than a random shift of any finite
configuration.

We remark that Lemma~\ref{L:evdom} is no longer true if the infection rates
fail to satisfy the first condition in (\ref{irr}). Indeed, if $A$ is defined
as in (\ref{irr}) and $\bigcup_{n\geq 0,\ m\geq 0}A^{-n}A^m\neq\La$, then we
can find sites $i,j\in\La$ such that there exists no site $k$ from which both
$i$ and $j$ can be infected. In particular, if we set $B:=\{i,j\}$, then
$\P[\exists i'\in\La\mbox{ s.t.\ }\eta^{\{0\}}_t\supset i'B]=0$ for all $t\geq
0$. For example, this happens if $\La$ is the free group with two generators,
say $g_1$ and $g_2$, $A=\{g_1,g_2\}$, and $B=\{g_1^{-1},g_2^{-1}\}$.\med

\detail{Note that $a(0,i)>0$ is equivalent to $a(k,ki)>0$ for any $k\in\La$.
  Set $C:=\bigcup_{n\geq 0}A^n$. Then a site $k$ can be infected from a site
  $j$ if and only if $k=ji$ (or, equivalently, $ki^{-1}=j$) for some $i\in
  C$. Now $\bigcup_{n\geq 0,\ m\geq 0}A^{-n}A^m\neq\La$ is equivalent to
  $C^{-1}C\neq\La$, which says that there exists a site $i$ that cannot be
  written as $i=j^{-1}k$ for some $j,k\in C$, which is equivalent to saying
  that there exists a site $i$ such that there do not exist sites $j,k\in C$
  such that $ik^{-1}=j^{-1}$, which is equivalent to saying that there exists
  a site $i$ such that there does not exists a site $k'$ such that both $i$
  and $0$ can be infected from $k'$ (i.e., such that $k'=ik^{-1}$ and
  $k'=0j^{-1}$ for some $j,k\in C$).}

\noi
{\bf Proof of Proposition~\ref{P:monhar}} Since the
$(\La,a,\de)$-contact process solves the martingale problem for $G$,
and $Gf=0$, the process $f(\eta^A_t)$ is a martingale. In particular:
\be
f(A)=\E[f(\eta^A_t)]\qquad(A\in\Pc_{\rm fin}(\La),\ t\geq 0).
\ee
Enumerate the elements of $\La$ an arbitrary way, and for
$A,B\in\Pc_{\rm fin}(\La)$, put
\be
\hat\imath_{A,B}:=\left\{\ba{ll}\min\{i\in\La:A\supset iB\}\quad
&\mbox{if }\{i\in\La:A\supset iB\}\mbox{ is nonempty,}\\
0&\mbox{otherwise.}\ea\right.
\ee
Since $f$ is monotone and shift invariant, we have, using Lemma~\ref{L:evdom},
\bc\label{xy}
\dis f(A)&=&\dis\lim_{t\to\infty}\E[f(\eta^A_t)]\\[5pt]
&\geq&\dis\limsup_{t\to\infty}\E[1_{\txt\{\exists i\in\La
\mbox{ s.t.\ }\eta^A_t\supset iB\}}
f(\hat\imath_{\eta^A_t,B}B)]\\[5pt]
&=&\dis f(B)\limsup_{t\to\infty}\P[\exists i\in\La
\mbox{ s.t.\ }\eta^A_t\supset iB]
=f(B)\rho(A)\qquad(A,B\in\Pc_{\rm fin}(\La)).
\ec
In particular, this shows that
\be
f(B)\leq\frac{f(\{0\})}{\rho(\{0\})}<\infty\qquad(B\in\Pc_{\rm fin}(\La)),
\ee
hence $f$ is bounded. Now let $A_n,B_m\in\Pc_{\rm fin}(\La)$ be
sequences such that $\rho(A_n)\to 1$ and $\rho(B_n)\to 1$. Then, by
(\ref{xy}),
\be
\liminf_{n\to\infty}f(A_n)\geq\liminf_{n\to\infty}f(B_m)\rho(A_n)
=f(B_m)\quad\forall m,
\ee
and therefore
\be\label{limsupinf}
\liminf_{n\to\infty}f(A_n)\geq\limsup_{m\to\infty}f(B_m).
\ee
This proves that the limit
\be
\lim_{\rho(A_n)\to 1}f(A_n)=:f(\infty)
\ee
exists and does not depend on the choice of the sequence $A_n$ with
$\rho(A_n)\to 1$. By the Markov property and continuity of the
conditional expectation with respect to increasing limits of
\si-fields (see Complement 10(b) from \cite[Section~29]{Loe63} or
\cite[Section~32]{Loe78}),
\be\label{loe}
\rho(\eta^A_t)
=\P\big[\eta^A_s\neq\emptyset\ \forall s\geq 0\,\big|\,\eta^A_t\big]\to
1_{\txt\{\eta^A_s\neq\emptyset\ \forall s\geq 0\}}\quad{\rm a.s.}
\quad\mbox{as }t\to\infty.
\ee
We conclude that, for all $A\in\Pc_{\rm fin}(\La)$, 
\be\ba{l}
\dis f(A)=\lim_{t\to\infty}\E[f(\eta^A_t)]\\[5pt]
\dis\quad=\P\big[\eta^A_t=\emptyset\mbox{ for some }t\geq 0\big]f(0)
+\P\big[\lim_{t\to\infty}\rho(\eta^A_t)=1\big]f(\infty)=
\rho(A)f(\infty),
\ec
which shows that $f$ is a scalar multiple of
$\rho$.\qed

\noi
The proof of Lemma~\ref{L:evdom} depends on two preparatory lemmas.
\bl{\bf(Local creation of finite configurations)}\label{L:loc}
For each $B\in\Pc_{\rm fin}(\La)$ and $t>0$, there exists a finite
$\De\sub\La$ and $j\in\La$ such that
\be\label{inDe}
\eps:=\P\big[\eta^{\{0\}}_t\supset jB\mbox{ and }
\eta^{\{0\}}_s\sub\De\ \forall 0\leq s\leq t\big]>0.
\ee
\el
{\bf Proof} It follows from assumption (\ref{irr}) that there
exists a site $j^{-1}\in\La$ with $\P\big[\eta^{\{j^{-1}\}}_t\supset B]>0$,
and therefore $\P\big[\eta^{\{0\}}_t\supset jB]>0$. Since
$\bigcup_{0\leq s\leq t}\eta^{\{0\}}_s$ is a.s.\ finite, we can choose
a finite but large enough $\De$ such that (\ref{inDe}) holds.\qed

\bl{\bf(Domination of finite configurations)}\label{L:dom}
For each $B\in\Pc_{\rm fin}(\La)$, $t>0$, and $A_n\in\Pc_{\rm fin}(\La)$
satisfying $\lim_{n\to\infty}|A_n|=\infty$, one has
\be
\lim_{n\to\infty}\P[\exists i\in\La\mbox{ s.t.\ }\eta^{A_n}_t\supset iB]=1.
\ee
\el
{\bf Proof} Let $\De$, $j$, and $\eps$ be as in Lemma~\ref{L:loc}. We can
find $\ti A_n\sub A_n$ such that $|\ti A_n|\to\infty$ as $n\to\infty$, and
for fixed $n$, the sets $(k\De)_{k\in\ti A_n}$ are disjoint. It follows that
\be\ba{l}
\P[\exists i\in\La\mbox{ s.t.\ }\eta^{A_n}_t\supset iB]\\[5pt]
\dis\qquad\geq1-\prod_{k\in\ti A_n}\big(1-\P\big[\eta^{\{k\}}_t\supset kjB
\mbox{ and }\eta^{\{k\}}_s\sub k\De\ \forall 0\leq s\leq t\big]\big)\\[5pt]
\dis\qquad=1-(1-\eps)^{|\ti A_n|}\asto{n}1,
\ec
where we have used (\ref{inDe}) and the fact that events concerning the
graphical representation in disjoint parts of space are independent.\qed

\noi
{\bf Proof of Lemma~\ref{L:evdom}} If $\de=0$, then obviously
$\lim_{t\to\infty}|\eta^A_t|=\infty$ a.s. If $\de>0$, then it is easy
to see that $\sup\{\rho(A):|A|\leq M\}<1$ for all $M<\infty$. Therefore,
by (\ref{loe}),
\be\label{exgro}
\eta^A_t=\emptyset\mbox{ for some }t\geq 0\quad\mbox{or}\quad|\eta^A_t|
\asto{t}\infty\qquad{\rm a.s.}
\ee
Fix $\emptyset\neq B\in\Pc_{\rm fin}(\La)$ and set
$\psi_t(A):=P[\exists i\in\La\mbox{ s.t.\ }\eta^A_t\supset iB]
\quad(A\in\Pc_{\rm fin}(\La),\ t\geq 0)$. Then, for each $t>0$,
\be
\lim_{T\to\infty}P[\exists i\in\La\mbox{ s.t.\ }\eta^A_T\supset iB]
=\lim_{T\to\infty}E[\psi_t(\eta^A_{T-t})]=\rho(A),
\ee
where we have used Lemma~\ref{L:dom} and (\ref{exgro}).\qed

\section{Proof of the main results}\label{S:proof}

\subsection{Exponentially growing processes}

In this section, we prove Theorem~\ref{T:main}~(d). Indeed, we prove
the following, more detailed result. Recall that by
Theorem~\ref{T:main}~(a) and Proposition~\ref{P:eigen}, there exists a
homogeneous eigenmeasure for the
$(\La,a^\dgg,\de)$-contact process, with eigenvalue
$r=r(\La,a^\dgg,\de)=r(\La,a,\de)$.
\bp{\bf(Exponential growth)}
Let $\La$ be a finite or countably infinite group, let
$a=(a(i,j))_{i,j\in\La}$ be infection rates satisfying (\ref{assum})
and let $\de\geq 0$. Let $\mu$ be any homogeneous
eigenmeasure of the $(\La,a^\dgg,\de)$-contact process,
with eigenvalue $r$, and let the function $v$ be defined in terms of
$\mu$ as in Lemma~\ref{L:eighar}. If $r>0$, then, for each
$A\in\Pc_{\rm fin}(\La)$, the limit
\be\label{WA}
W_A:=\lim_{t\to\infty}e^{-rt}v(\eta^A_t)
\ee
exists a.s., and satisfies $\E[W_A]=v(A)$. If the infection rates satisfy
(\ref{irr}), then moreover
\be
\P[W_A>0]=\P[\eta^A_t\neq\emptyset\ \forall t\geq 0].
\ee
\ep
{\bf Proof} Our proof follows a strategy that is familiar from the
theory of supercritical branching processes. Using duality (see
(\ref{dual})) and the fact that $\mu$ is an eigenmeasure, we see that
\be\ba{l}\label{Stv}
\dis\E[v(\eta^A_t)]=\int\mu(\di B)\P[\eta^A_t\cap B\neq\emptyset]
=\int\mu(\di B)\P[A\cap\eta^{\dgg\,B}_t\neq\emptyset]\\[5pt]
\dis\quad=e^{rt}\int\mu(\di B)\P[A\cap B\neq\emptyset]=e^{rt}v(A),
\ec
so, by the Markov property of $(\eta^A_t)_{t\geq 0}$, the process
$(e^{-rt}v(\eta^A_t))_{t\geq 0}$ is a martingale. Every nonnegative
martingale converges, so for each $A\in\Pc_{\rm fin}(\La)$, there
exists a random variable $W_A$ such that (\ref{WA}) holds.

To prove that $\E[W_A]=v(A)$, it suffices to show that the random
variables $\{e^{-rt}v(\eta^A_t):t\geq 0\}$ are uniformly
integrable. By Proposition~\ref{P:cov}, the variance of $v(\eta^A_t)$
is given by
\be
\var(v(\eta^A_t))=2\int_0^t\E\big[\Ga(S_sv,S_sv)(\eta^A_{t-s})\big]\di s,
\ee
where $S_t$ and $\Ga$ are defined in (\ref{St}) and (\ref{Gadef}).
Formula (\ref{Stv}) tells us that $S_tv=e^{rt}v$, so
\be\label{varv}
\var(v(\eta^A_t))=2\int_0^t\E\big[\Ga(e^{rs}v,e^{rs}v)(\eta^A_{t-s})\big]\di s
=2\int_0^t\E\big[\Ga(v,v)(\eta^A_{t-s})\big]e^{2rs}\di s.
\ee
It is not hard to see that for any $f,g\in\Si(\Pc_{\rm fin}(\La))$
and $A\in\Pc_{\rm fin}(\La)$,
\bc
\dis\Ga(f,g)(A)
&=&\dis\ffrac{1}{2}\sum_{ij}a(i,j)1_{\{i\in A\}}1_{\{j\not\in A\}}
\big(f(A\cup\{j\})-f(A)\big)\big(g(A\cup\{j\})-g(A)\big)\\[5pt]
&&\dis+\ffrac{1}{2}\de\sum_i1_{\{i\in A\}}
\big(f(A\beh\{i\})-f(A)\big)\big(g(A\beh\{i\})-g(A)\big).
\ec
Without loss of generality we can normalize $v$ such that
$v(\{0\})=1$. Then, by monotonicity and subadditivity (see
Lemma~\ref{L:eighar}), $0\leq v(A\cup\{j\})-v(A)\leq 1$ for all $j,A$,
and therefore
\be
\Ga(v,v)(A)\leq\ffrac{1}{2}\sum_{ij}a(i,j)1_{\{i\in A\}}1_{\{j\not\in A\}}
+\ffrac{1}{2}\de\sum_i1_{\{i\in A\}}\leq\ffrac{1}{2}\big(|a|+\de\big)|A|.
\ee
Inserting this into (\ref{varv}) yields
\bc
\dis\var(e^{-rt}v(\eta^A_t))
&\leq&\dis\big(|a|+\de\big)e^{-2rt}
\int_0^t\E\big[|\eta^A_{t-s}|\big]e^{2rs}\di s\\[5pt]
&\leq&\dis\big(|a|+\de\big)
\int_0^t\E\big[|\eta^A_{t-s}|\big]e^{-2r(t-s)}\di s.
\ec
Since $r$ is the exponential growth rate of $\eta^A$ and $r>0$, we can
find $K<\infty$ such that $\E\big[|\eta^A_t|\big]\leq
Ke^{\frac{3}{2}rt}$ $(t\geq 0)$. It follows that
\be
\var(e^{-rt}v(\eta^A_t))
\leq(|a|+\de)K\!\int_0^\infty e^{-\frac{1}{2}rs}\di s<\infty\qquad(t\geq 0),
\ee
which proves the required uniform integrability.

Set $f(A):=\P[W_A>0]$ and recall that $\rho(A):=\P[\eta^A_t\neq\emptyset\
\forall t\geq 0]$. Obviously $f\leq\rho$. We have just shown that
$\P[W_A>0]>0$ if $A\neq\emptyset$, so $\rho(A)\geq f(A)>0$ for each
$A\neq\emptyset$. Assuming that the infection rates satisfy
(\ref{irr}), we claim that $f=\rho$. We observe that
\be\label{loef}
f(\eta^A_t)=\P\big[\lim_{s\to\infty}e^{-rs}\eta^A_s>0\,\big|\,\eta^A_t\big].
\ee
In particular, this shows that $(f(\eta^A_t))_{t\geq 0}$ is a
martingale, hence $Gf=0$. It is easy to see that $f$ is shift-invariant,
monotone, bounded, and satisfies $f(\emptyset)=0$, so applying
Proposition~\ref{P:monhar}, we see that $f=c\rho$ for some $c\geq 0$.
Since $f\leq\rho$, we have $c\leq 1$.

By continuity of the conditional expectation with respect to
increasing limits of \si-fields (compare (\ref{loe})), the right-hand
side of (\ref{loef}) converges a.s.\ to the indicator function of the
event that $W_A>0$. Since this event has positive probability, the
event $\lim_{t\to\infty}f(\eta^A_t)=1$ has positive probability. In
particular, this shows that for each $\eps>0$ there exists a finite
set $B$ with $f(B)\geq 1-\eps$. This is possible only if the
constant $c$ in the equation $f=c\rho$ satisfies $c\geq 1$.\qed

\subsection{Subexponential lattices}

\newcommand{\lc}{h}

{\bf Proof of Theorem~\ref{T:main}~(e)} Consider a branching
process on $\La$, started with one particle in the origin, where a
particle at $i$ produces a new particle at $j$ with rate $a(i,j)$, and
each particle dies with rate $\de$. Let $B_t(i)$ denote the number of
particles at site $i\in\La$ and time $t\geq 0$. It is not hard to see
that $\eta^{\{0\}}$ and $B$ may be coupled such that
\be\label{brest}
1_{\eta^{\{0\}}_t}\leq B_t\qquad(t\geq 0).
\ee
Let $(\xi_t)_{t\geq 0}$ be a random walk on $\La$ that jumps from $i$
to $j$ with rate $a(i,j)$, started in $\xi_0=0$.  Then it is not hard
to see that (compare \cite[Proposition~I.1.21]{Lig99})
\be\label{walk}
\E[B_t(i)]=\P[\xi_t=i]e^{(|a|-\de)t}\qquad(i\in\La,\ t\geq 0).
\ee
\detail{Indeed, if $b_t(i)=E[B_t(i)]$, then
\[
\dif{t}b_t(i)=\sum_ja(j,i)b_t(j)-b_t(i)
=(|a|-\de)b_t(i)+\sum_ja(j,i)b_t(j)-|a|b_t(i).
\]}
Let $\lc>0$ be a constant, to be determined later. It follows from
(\ref{brest}) and (\ref{walk}) that
\bc\label{split}
\dis\E\big[|\eta^{\{0\}}_t|\big]&\leq&
\dis\sum_i\big(1\wedge\P[\xi_t=i]e^{(|a|-\de)t}\big)\\[5pt]
&=&|\{i\in\La:|i|\leq\lc t\}|+\P[|\xi_t|>\lc t]e^{(|a|-\de)t}\qquad(t\geq 0).
\ec
Let $(Y_i)_{i\geq 1}$ be i.i.d.\ $\N$-valued random variables with
$\P[Y_i=k]=\frac{1}{|a|}\sum_{j:\ |j|=k}a(0,j)$ $(k\geq 0)$, let $N$ be
a Poisson-distributed random variable with mean $|a|$, independent of
the $(Y_i)_{i\geq 1}$, and let $(X_m)_{m\geq 1}$ be i.i.d.\ random
variables with law
$\P[X_m\in\cdot\,]=\P[\sum_{i=1}^NY_i\in\cdot\,]$. Since the random walk
$\xi$ makes jumps whose sizes are distributed in the same way as the
$Y_i$, and the number of jumps per unit of time is Poisson distributed
with mean $|a|$, it follows that
\be\label{linsplit}
\P[|\xi_t|>\lc t]\leq\P\Big[\frac{1}{\lceil t\rceil}
\sum_{m=1}^{\lceil t\rceil}X_m>\lc\frac{t}{\lceil t\rceil}\Big]\qquad(t>0),
\ee
where $\lceil t\rceil$ denotes $t$ rounded up to the next integer.
By our assumptions,
\be
\E\big[\ex{\eps X_m}\big]=\E\big[\ex{\eps\sum_{i=1}^NY_k}\big]
=e^{-|a|}\sum_{n=0}^\infty\frac{|a|^n}{n!}\E\big[\ex{\eps Y_1}\big]^n
=\ex{-|a|(1-\E[e^{\eps Y_1}])}<\infty,
\ee
for some $\eps>0$. Therefore, by \cite[Theorem~2.2.3 and Lemma~2.2.20]{DZ98},
for each $R>0$ there exists a $\lc>0$ and $K<\infty$ such that
\be\label{Rex}
\P\Big[\frac{1}{n}\sum_{m=1}^{n}X_m>\lc\Big]\leq K\ex{-nR}\qquad(n\geq 1).
\ee
Choosing $\lc$ such that (\ref{Rex}) holds for some $R>|a|-\de$ yields,
by (\ref{linsplit})
\be
\lim_{t\to\infty}\P\big[|\xi_t|>\lc t\big]e^{(|a|-\de)t}=0.
\ee
Inserting this into (\ref{split}) we find that the exponential growth
rate $r=r(\La,a,\de)$ satisfies
\be
r\leq\limsup_{t\to\infty}\ffrac{1}{t}\log|\{i\in\La:|i|\leq\lc t\}|=0,
\ee
where we have used that the group $\La$ has subexponential growth.\qed

\subsection{Nonamenable lattices}\label{S:Palm}

In this section, we prove Theorem~\ref{T:main}~(f) and
Corollary~\ref{C:crit}. We start by introducing some notation.
If $R$ is a nonnegative real random variable, defined on some
probability space $(\om,\Fi,\P)$, and $0<\E[R]<\infty$, then
we define the {\em size-biased law} $\ov\P_R$ associated with $R$ by
\be\label{size}
\ov\P_R(\Ai):=\frac{\E[1_\Ai R]}{\E[R]}\qquad(\Ai\in\Fi).
\ee
If $\De$ is a $\Pc_{\rm fin}(\La)$-valued random variable, defined on
some probability space $(\om,\Fi,\P)$, such that
$0<\E[|\De|]<\infty$, then we define a probability law
$\hat\P=\hat\P_\De$ on the product space $\om\times\La$ by
\be\label{Camp}
\hat\P_\De(\Ai\times\{i\}):=\frac{\P(\{i\in\De\}\cap\Ai)}{\E[|\De|]}
\qquad(\Ai\in\Fi,\ i\in\La).
\ee
We call $\hat\P_\De$ the {\em Campbell law} associated with $\De$. It
is not hard to see that the projection of $\P_\De$ onto $\om$ is the
size-biased law $\ov\P_{|\De|}$. Moreover, if we let $\iota(\oo,i):=i$
denote the projection from $\om\times\La$ to $\La$ and we use the
symbol $\De$ to denote (also) the random variable on $\om\times\La$
defined by $\De(\oo,i):=\De(\oo)$, then 
\be
\hat\P_\De\big[\iota=i\,\big|\,\De\big]=\ffrac{1}{|\De|}1_\De(i),
\ee
i.e., conditional on $\De$, the site $\iota$ is chosen with equal
probabilities from all sites in $\De$. We may view $\iota$ as a
`typical' element of~$\De$. Campbell laws are closely related to the
more widely known Palm laws; both play an important role in the theory
of branching processes (see, e.g., \cite{Win99}).

The next lemma relates Campbell laws to things we have been
considering so far. Note that if $\mu$ is a locally finite measure on
$\Pc_+(\La)$ and $\int1_{\{0\in A\}}\mu(\di A)>0$, then the
conditional law
\be
\mu(\di A\,|\,0\in A)
:=\frac{1_{\{0\in A\}}\mu(\di A)}{\int1_{\{0\in B\}}\mu(\di B)}
\ee
is a well-defined probability law.
\bl{\bf(Campbell law)}\label{L:Camp}
Let $\eta$ be a $(\La,a,\de)$-contact process. For each $t\geq 0$, let
$\mu_t$ be defined as in (\ref{mut}) and let
$\hat\P_t:=\hat\P_{\eta^{\{0\}}_t}$ be the Campbell law associated
with $\eta^{\{0\}}_t$. Then
\be
\mu_t(\di A\,|\,0\in A)=\hat\P_t\big[\iota^{-1}\eta^{\{0\}}_t\in\di A\big]
\ee
\el
{\bf Proof} This follows by writing
\be\ba{l}
\dis\hat\P_t\big[\iota^{-1}\eta^{\{0\}}_t\in\di A\big]
=\sum_i\hat\P_t\big[i^{-1}\eta^{\{0\}}_t\in\di A,\ \iota=i\big]\\[5pt]
\dis\quad
=\sum_i\frac{\P\big[i^{-1}\eta^{\{0\}}_t\in\di A,\ i\in\eta^{\{0\}}_t\big]}
{\E[|\eta^{\{0\}}_t|]}
=\frac{\sum_i\P\big[i^{-1}\eta^{\{0\}}_t\in\di A,\ i\in\eta^{\{0\}}_t\big]}
{\sum_i\P[i\in\eta^{\{0\}}_t]}\\[10pt]
\dis\quad
=\frac{\sum_i\P\big[i^{-1}\eta^{\{0\}}_t\in\di A,
\ 0\in i^{-1}\eta^{\{0\}}_t\big]}
{\sum_i\P[0\in i^{-1}\eta^{\{0\}}_t]}
=\frac{\sum_i\P\big[\eta^{\{i^{-1}\}}_t\in\di A,
\ 0\in\eta^{\{i^{-1}\}}_t\big]}
{\sum_i\P[0\in\eta^{\{i^{-1}\}}_t]}\\[10pt]
\dis\quad
=\frac{\sum_j\P\big[\eta^{\{j\}}_t\in\di A,\ 0\in\eta^{\{j\}}_t\big]}
{\sum_j\P[0\in\eta^{\{j\}}_t]}
=\frac{1_{\{0\in A\}}\sum_j\P\big[\eta^{\{j\}}_t\in\di A\big]}
{\int1_{\{0\in B\}}\sum_j\P\big[\eta^{\{j\}}_t\in\di B\big]}\\[10pt]
\dis\quad
=\frac{1_{\{0\in A\}}\mu_t(\di A)}{\int1_{\{0\in B\}}\mu_t(\di B)}
=\mu_t(\di A\,|\,0\in A).
\ec
\qed

\noi
The next proposition is a direct consequence of Theorem~\ref{T:zero} and
Corollary~\ref{C:conv}. Note that by the remark below Corollary~\ref{C:conv},
we expect the convergence in (\ref{Campconv}) to hold also when the
$\tau_\ga$ are replaced by deterministic times tending to infinity.
\bp{\bf(Convergence of Campbell laws)}\label{P:Campconv}
Assume that the $(\La,a,\de)$-contact process has a nontrivial upper
invariant measure $\ov\nu$, that its exponential growth rate
$r(\La,a,\de)$ is zero, and that the infection rates satisfy
(\ref{irr}). Let $\eta^{\{0\}}$ be the $(\La,a,\de)$-contact process
started in $\{0\}$ and for $\ga\geq 0$, let $\tau_\ga$ be an
exponentially distributed random variable with mean $\ga$, independent
of $\eta^{\{0\}}$. For each $\ga\geq 0$, let $\hat\P_\ga=\hat
\P_{\eta^{\{0\}}_{\tau_\ga}}$ be the Campbell law associated with
$\eta^{\{0\}}_{\tau_\ga}$. Then
\be\label{Campconv}
\hat\P_\ga\big[\iota^{-1}\eta^{\{0\}}_{\tau_\ga}\in\di A\big]
\Asto{\ga}\ov\nu(\di A\,|\,0\in A),
\ee
where $\Rightarrow$ denotes weak convergence of probability measures.
\ep
{\bf Proof} For $\la>0$, let $\hat\mu_\la$ denote the Laplace transform
of $\mu_t$, defined in (\ref{Lap}). In analogy with Lemma~\ref{L:Camp}, it
is straightforward to check that
\be\label{condlaw}
\hat\mu_\la(\di A\,|\,0\in A)=
\hat\P_{1/\la}\big[\iota^{-1}\eta^{\{0\}}_{\tau_{1/\la}}\in\di A\big]
\qquad(\la>0).
\ee
By Theorem~\ref{T:zero} and Corollary~\ref{C:conv}, the measures
$\hat\mu_\la$, suitably rescaled, converge vaguely to $\ov\nu$ as
$\la\down 0$. By (\ref{condlaw}), this implies the weak convergence in
(\ref{Campconv}).\qed

\noi
The next proposition shows that if the assumptions of
Proposition~\ref{P:Campconv} are satisfied, then $\La$ must be
amenable.
\bp{\bf(Campbell laws and amenability)}\label{P:Camen}
Let $\La$ be a countable group, let $B_n$ be random nonempty, finite subsets
of $\La$, and conditional on $B_n$, let $\iota_n$ be chosen with equal
probabilities from the sites in $B_n$. Let $B$ be a random subset of $\La$
whose law is nontrivial and homogeneous. Assume that
\be\label{Camen}
\P\big[\iota_n^{-1}B_n\in\cdot\,\big]\Asto{n}\P[B\in\cdot\,\big|\,0\in B].
\ee
Then $\La$ must be amenable.
\ep
{\bf Proof} The idea behind Proposition~\ref{P:Camen} is easy to explain:
formula (\ref{Camen}) says that for large $n$, the set $B_n$ looks like a
random finite piece cut out of the spatially homogeneous configuration $B$,
such that most points in this piece are far from the boundary. This
contradicts nonamenability, since in any finite subset of a nonamenable group,
a positive fraction of the points must lie near the boundary.

To make this idea rigorous, we proceed as follows. Assume that $\La$ is
nonamenable. Then there exists a finite nonempty $\De\sub\La$ and $\eps>0$
such that $|(A\De)\symdif A|\geq\eps|A|$ for all finite nonempty
$A\sub\La$. Without loss of generality we may assume that $\De$
is symmetric. Let $(\xi_m)_{m\geq 0}$ be a random walk in $\La$, independent
of $B_n$ and $B$, starting in $\xi_0=0$, that jumps from a point $i$ to a
point $ij$ with probability $|\De|^{-1}1_{\{j\in\De\}}$. Then (\ref{Camen})
implies that
\be\label{limxi}
\P[\iota_n\xi_m\in B_n]\Asto{n}\P[\xi_m\in B\,|\,\xi_0\in B]\qquad(m\geq 0).
\ee
By the stationarity of the process $(1_{\{\xi_m\in B\}})_{m\geq 0}$, one has
\be\label{nonzero}
\limsup_{m\to\infty}\P[\xi_m\in B\,|\,\xi_0\in B]>0.
\ee
On the other hand, we will show that the nonamenability of $\La$ implies that
\be\label{unizero}
\lim_{m\to\infty}\sup_{n\geq 0}\P[\iota_n\xi_m\in B_n]=0,
\ee
which with (\ref{nonzero}) leads to a contradiction in (\ref{limxi}). Let
$\ell^2(\La)$ be the Hilbert space of square summable real functions on $\La$,
equipped with the inner product $\li x,y\re:=\sum_ix(i)y(i)$, let
$P^n(i,j):=\P[\xi_n=j\,|\,\xi_0=i]$ and $P^nx(i):=\sum_jP^n(i,j)x(j)$. Then,
by the fact that nearest-neighbor random walk on any nonamenable Cayley graph
has a spectral gap (see \cite{Kes59} or \cite[Thm~6.7]{LP08}), there exists a
$0<\tet<1$ such that
\be
|B_n|\,\P[\iota_n\xi_m\in B_n]=\sum_{i\in B_n}\sum_{j\in B_n}P^m(i,j)
=\li1_{B_n}, P^m1_{B_n}\re\leq\tet^m\li1_{B_n},1_{B_n}\re=\tet^m|B_n|,
\ee
which proves (\ref{unizero}).\qed

\detail{If $0\not\in\De$, then set $\ti\De:=\De\cup\{0\}$ and
observe that $(A\ti\De)\symdif A=(A\De)\beh A$ while $|A\beh(A\De)|=|\{i\in
A:ij\not\in A\ \forall j\in\De\}|\leq|\{i\in A:\exists j\in\De\mbox{
  s.t.\ }ij\not\in A\}|=|\{ij\not\in A:i\in A,\ j\in\De\}|=|(A\De)\beh A|$,
hence setting $\ti\eps:=\eps/2$ yields $|(A\ti\De)\symdif A|\geq\ti\eps|A|$
for all finite nonempty $A\sub\La$. If $0\in\De$ but $\De$ is not symmetric
then setting $\ti\De:=\De\cup\De^{-1}$ yields $|A\ti\De\beh A|\geq|A\De\beh
A|\geq\eps|A|$ for each finite nonempty $A\sub\La$.}

\noi
{\bf Proof of Theorem~\ref{T:main}~(f)} Assume (\ref{irr}). Assume
that the $(\La,a,\de)$-contact process survives and that its
exponential growth rate $r(\La,a,\de)$ is zero. Then the
$(\La,a^\dgg,\de)$-contact process has a nontrivial upper invariant
law, and, by Theorem~\ref{T:main}~(a), $r(\La,a^\dgg,\de)=0$.
Therefore, by Propositions~\ref{P:Campconv} and \ref{P:Camen}, $\La$
must be amenable.\qed

\noi
{\bf Proof of Corollary~\ref{C:crit}} Let $\Si:=\{\de\geq 0:\mbox{the
$(\La,a,\de)$-contact process survives}\}$. Note that $\Si$ is
nonempty since $0\in\Si$. By Theorem~\ref{T:main}~(d) and (f),
$\Si=\{\de\geq 0:r(\La,a,\de)>0\}$. By Theorem~\ref{T:main}~(b), the
function $\de\to r(\La,a,\de)$ is continuous, hence $\Si$ is an open
subset of $\half$. Hence, by monotonicity, $\Si=[0,\de_{\rm c})$, where
$\de_{\rm c}:=\sup\{\de\geq 0:$ the $(\La,a,\de)$-contact process
survives$\}$ satisfies $\de_{\rm c}>0$.\qed

\end{document}